\documentclass[12pt]{article}
\makeatletter \oddsidemargin  -.1in \evensidemargin -.1in
\newcommand{\singlespacing}{\let\CS=\@currsize\renewcommand{\baselinestretch}{1}\tiny\CS}
\newcommand{\oneandahalfspacing}{\let\CS=\@currsize\renewcommand{\baselinestretch}{1.25}\tiny\CS}
\newcommand{\doublespacing}{\let\CS=\@currsize\renewcommand{\baselinestretch}{1.35}\tiny\CS}

\newtheorem{rule-def}[theorem]{Rule}

\RequirePackage[dvips]{graphicx} \textheight 22.5cm
\usepackage{latexsym,epsfig,enumerate,amsmath,amsfonts,amssymb,amsbsy,amsopn,mathrsfs}
\usepackage{subfigure}
\usepackage{graphicx}

\begin{document}

\title{\bf Results of Brocard-Ramanujan problem on diophantine
  equation $n!+1=m^2$} \author{{\small Somnath Maiti$^{1,2}$}
  \thanks{Sections 1-2.6 carried out by April 2020 at The LNM Institute of Information Technology, Jaipur, India}\\\it $^{1}$ Department of Mathematics, Rajendra College\\\it (A Constituent Unit of Jai Prakash University), Chapra 841301, Bihar, India\\Email address: \it maiti0000000somnath@gmail.com\\ \it
  $^{2}$Department of Mathematical Sciences, Indian Institute of
  Technology (BHU),\\ Varanasi 221005, India} \date{}
\maketitle \noindent \doublespacing
\vspace{-0.5cm}
\begin{abstract}
The Brocard-Ramanujan problem pertaining to the diophantine equation
$n!+1=m^2$, a famously unsolved problem, deals with finding the
integer solutions to the equation. Nobody has discovered any new solution of the problem beyond $n=4,~5$ and $7$ although many of us have tried it. Bruce Berndt and William Galway \cite{Berndt} had not found any new solution in 2000 by extensive computer search for a solution with $n$ up to $10^9$. The purpose of this study is to show that the solutions should satisfy some necessary and/or sufficient conditions. If $\sqrt{n!}=k+\epsilon,~n>1,~0<\epsilon<1$; then it has solution if and only if $n!=k(k+2)$ and $\epsilon,~k$  are strictly monotonic increasing. It has only finitely many solutions which is not based on any conjecture or previous research on the Brocard-Ramanujan problem. For the new solution of Brocard-
Ramanujan problem ($n\ge 10^5$), the value of $\epsilon$ should be more than $0.999 \cdots 905915$ (digit 0 is coming after 228287 numbers of 9 digit, which takes more than 66 pages in (LibreOffice Writer) indicating almost impossibility of new solution. If we consider $n\geq 10^9$, I am unable to calculate the said numbers of 9 digit in the value of $\epsilon$
in my personal laptop (with 8GB Ram) using MATHEMATICA 8.  Finally, it has been claimed to discover that the problem has no further solution.  
\end{abstract}

\section{Introduction}
The Brocard-Ramanujan problem was reported by Henri Brocard
\cite{Brocard1,Brocard2} in the two written articles in 1876 and 1885
asking to search positive integers $n$ and $m$ satisfying the equation
$n!+1=m^2$, where $n!$ is the factorial. Later in 1913, unaware of
Brocard's query, S. Ramanujan \cite{Ramanujan} reported the same
problem in the following composition: ``The number $n!+1$ is a perfect
square for the values 4, 5, 7 of $n$. Find other values.''

It is believed that the only solutions of the problem are $(m,n)=(5,
4), (11, 5), (71, 7)$. The solutions of the problem in the pair form
$(m,n)$ is known as Brown's numbers. As reported by
A. G$\acute{e}$rardin \cite{Gerardin} with arguments in 1906, the
problem has no solutions for $7< n<25$ and $m$ must have at least 20
digits if $m>71$ for the solution. The computations of $n!$ up to
$n=63$ carried out by Hansraj Gupta \cite{Gupta} in 1935 provided no
further solutions of the problem. Considering more general diophantine
equations $x^k\pm y^k=n!$ and $m!\pm n!=x^k$, P. Erd\H{o}s and
R. Obl\'{a}th \cite{Erdos} reported that most of these have only
finitely many solutions. For example, the equation $x^k-1=n!$ has no
solutions if $k>1$ and possible exception of $k=2$ and $k=4$. They
proved that the number of solutions is finite for the case of $k=4$
i.e. for the equation $x^4-1=n!$. Pollack and Shapiro \cite{Pollack}
finally showed that the equation $x^4-1=n!$ has actually no solutions.

In 1993, Marius Overholt \cite{Overholt} proved that ``a weak form of
the unproved diophantine inequality known as Szpiro's conjecture would
imply that there are only finitely many solutions'' of the
problem. The weak form of Szpiro's conjecture, a special case of the
ABC conjecture, can be defined as follows: Let
$N_0(n)=\displaystyle\prod_{p\mid n}p$, in which $p$ indicates a
prime; $a$, $b$, and $c$ denote positive integers together with
relatively prime in pairs and follow the equation $a+b=c$. Then the
weak form of Szpiro's conjecture implies that there exists a constant
$k$ such that $|abc|\le [N_0(abc)]^k$.

In 1996, A. Dabrowski \cite{Dabrowski} carried out the
extension of the result, under the identical condition, to all
equations of the type $n!+A=m^2$ in terms of an arbitrary natural
number $A$. Under the identical condition or even a harder constitution, the
result of Overholt was further carried forward by several researchers
\cite{Luca,Gica1,Gica2} to equations of the form $n!=P(m)$ or even $n(n-1)\cdots
(n-k+1)=P(m)$, in which $P$ is a polynomial of degree $\geq$2 with
integer coefficients, and $k,~m,~n$ are integers with $1\leq k\leq n$.

In 2000, Bruce Berndt and William Galway \cite{Berndt} found no solution of the
problem based on computations up to $n=10^9$, reported a relation
between Legendre’s symbol $(\frac{a}{p})$ and the solutions of
$n!+1=m^2$ indicating $(\frac{n!+1}{p})=1~\text{or}~0$.

On the basis of consideration of Szpiro's conjecture, Dufour and Kihel
\cite{Dufour} in 2004 proved that Hall's conjecture weak form implies
that $n!+A=m^2$ (if the integer $A$ is not a $q$-th power of an
integer) can only have finitely many solutions, where Hall's conjecture
defines a special form of the abc-conjecture. 

In 2008, Kihel and Luca \cite{Kihel} presented some variants of the
Brocard-Ramanujan problem considering the form $\displaystyle x^p\pm
y^p=\prod_{k=1~(k\nmid n)}^{n} k$ and showed that there are only
finitely many integer solutions $(x,~y,~p,~n)$ if $p~(\ge 3)$ a prime
number and gcd$(x,y)=1$. 

In 2013, Liu \cite{Liu} wrote a Master of Science thesis on Brocard’s
problem and its variations in Waikato University on the basis of
reviews for the old methods in order to solve various forms of the
Brocard’s problem and presented some proofs in his own way. Concerning
variations of Brocard’s problem, Maciej Ulas \cite{Ulas} considered
the equation $n!+A=m^2$ and proved that there were infinitely many
$A’$s such that the equation $n!+A=m^2$ has at least three solutions
in natural number. Dabrowski and Ulas \cite{Dabrowski2} investigated
Diophantine equations of the form $m^2=px_n+q$, in which $p$, $q$ are
fixed integers, $x_n=g(1)g(2)\cdots g(n)$ and $g:
\mathbb{N}_+\rightarrow \mathbb{N}_+$ is an increasing function. They
proved that the equations had at least four solutions in positive
integers $m$ and $n$.

Maximum number of the results of the Brocard-Ramanujan problem and its 
variations has been derived by using the ABC conjecture or its
variants \cite{Dabrowski2}. A question is arising naturally that how the results can be proved without the help of this assumption \cite{Dabrowski2}? In view of the above, one of the
purpose of this study is to show the finiteness of the number of
solutions of $n!+1=m^2$, which is not based on any conjecture or
previous study and the solutions should satisfy some necessary and/or sufficient conditions. The proof of the existence of new solution of equation (\ref{brocard_ramanujan_problem}) seems to be as distant as the establishment of the nonexistence of odd perfect numbers and the two problems although not equivalent are not dissimilar. 

We are tempted to believe the stronger statement that
Brocard–Ramanujan problem has no other solutions except $(m,n)=(5, 4),~(11, 5),~(71, 7)$; however, it can be claimed that this study have been able to establish this stronger statement. Moreover, results on Brocard–Ramanujan problem show that if $\sqrt{n!}=k+\epsilon,~n>1,~0<\epsilon<1$; then it has solution if and only if $n!=k(k+2)$ and $\epsilon,~k$  are strictly monotonic increasing. For new solution of Brocard-
Ramanujan problem considering $n\ge 10^5$, the value of $\epsilon$ should be more than $0.999 \cdots 905915$ (digit 0 is coming after 228287 numbers of 9 digit, which takes more than 66 pages in LibreOffice Writer) concluding the impossibility of new solution. I am unable to calculate the numbers of said 9 digit in the value of $\epsilon$
in my personal laptop (with 8GB Ram) using MATHEMATICA 8 when $n\geq 10^9$.

\section{Results and Discussion}
\subsection{Theorem}
\label{brocard_ramanujan_1st_theorem}
The natural number $m$ will be a solution of the diophantine equation 
\begin{equation}
n!+1=m^2;~n,m\in \mathbb{N}
\label{brocard_ramanujan_problem}
\end{equation}
if and only if $m=[\sqrt{n!}]+1$, where for $x\in\mathbb{R},~[x]$
indicates greatest integer function. Moreover for the solution of
equation (\ref{brocard_ramanujan_problem}), we can express $n!$ as
$n!=2^{\beta_1}\times \text{odd number}=(2 a)\times (2^{\beta_1-1}b)$ with $a$ and
$b$ relatively odd natural numbers satisfying gcd$(a,2^{\beta_1-1}b)=1$ and
the difference of $2a$ and $2^{\beta_1-1}b$ will be 2 with $\beta_1=[\frac{n}{2}]+[\frac{n}{4}]+[\frac{n}{8}]+\cdots$.

{\bf Proof}: The equation (\ref{brocard_ramanujan_problem}) can be
written as $n!=m^2-1=(m-1)(m+1)$.

Then we can express it in the form:
\begin{equation}
n!=[\sqrt{n!}]([\sqrt{n!}]+2)~~~~~~~~~~~~~~~~~~~~~~~~~~~~~~~~~~~~~~~~~~~~~~~~~~~~~~~~~~~~~~
\label{brocard_ramanujan_problem_solution_form}
\end{equation} 

Or, $[\sqrt{n!}]^2+2[\sqrt{n!}]-n!=0$. 

Or, $[\sqrt{n!}]=\frac{-2\pm\sqrt{4+4n!}}{2}=-1\pm \sqrt{1+n!}$. 

So, $[\sqrt{n!}]=\sqrt{1+n!}-1$. 
\begin{equation}
\text{Or,}~[\sqrt{n!}]+1=\sqrt{1+n!}=\sqrt{m^2}=m~~~~~~~~~~~~~~~~~~~~~~~~~~~~~~~~~~~~~~~~~~~~~~~~~~~~~~~~~~
\label{brocard_ramanujan_problem_expression_m}
\end{equation} 
Thus if $m$ is a solution of equation
(\ref{brocard_ramanujan_problem}), then it will be in the form
$m=[\sqrt{n!}]+1$. Hence the condition is necessary.

Conversely, let $m=[\sqrt{n!}]+1$. Then $(m-1)(m+1)=[\sqrt{n!}]([\sqrt{n!}]+2)$.

Or, $n!=m^2-1=[\sqrt{n!}]([\sqrt{n!}]+2)$ i.e. $n!+1=([\sqrt{n!}]+1)^2$.

Then, Brown's number $([\sqrt{n!}]+1,n)$ will be a solution of the equation
(\ref{brocard_ramanujan_problem}). Thus, the condition is sufficient.

Again from equations (\ref{brocard_ramanujan_problem_solution_form})
and (\ref{brocard_ramanujan_problem_expression_m}), we get
$(m-1)(m+1)=[\sqrt{n!}]([\sqrt{n!}]+2)=n!=2^e\times \text{odd number}$
for some $e\in \mathbb{N}$ and gcd$(m-1,m+1)=d$ where
$d=1~\text{or}~2$ since $d\mid (m-1)$ and $d\mid (m+1)$ implies that
$d\mid \{(m+1)-(m-1)\}$ i.e. $d\mid 2$. If $d=1$, then one of
$(m-1),~(m+1)$ and hence other also will be odd natural number. Thus
$n!=\text{odd number}~\text{for}~n>1$, which is a absurd result. Hence
$d=2$ and we can express $n!$ as $n!=2^k\times \text{odd number}=(2
a)\times (2^{e-1}b)$ with $a$ and $b$ relatively odd natural numbers
satisfying gcd$(a,2^{k-1}b)=1$ and the difference between $2a$ and
$2^{k-1}b$ will be 2.

\subsubsection{Note}
(i) If we consider the only remaining option for the expression of $n!$ as
$n!=([\sqrt{n!}]-1)([\sqrt{n!}]+1)$, then $[\sqrt{n!}]^2-1=n!$.

Or, $[\sqrt{n!}]=\sqrt{1+n!}$.
\begin{equation}
\text{Thus,}~~~~[\sqrt{n!}]=\sqrt{1+n!}=\sqrt{m^2}=m.~~~~~~~~~~~~~~~~~~~~~~~~~~~~~~~~~~~~~~~~~~~~~~~~~~~~~~~~~~~~
\label{brocard_ramanujan_n!_another_expression}
\end{equation}
 Now if $\sqrt{n!}=r\in\mathbb{N}$, then $n!=r^2$ and hence
 $n!+1=r^2+1\neq m^2,~m\in \mathbb{N}$. 

Thus $\sqrt{1+n!}=m$ implies that $[\sqrt{n!}]<m$, a contradiction of
the equation (\ref{brocard_ramanujan_n!_another_expression}). Hence
the expression $n!=([\sqrt{n!}]-1)([\sqrt{n!}]+1)$ can not give any
solution of $n!+1=m^2;~n,m\in \mathbb{N}$.

(ii) If $k=m-1$, then $k^2<n!=k(k+2)<(k+1)^2$. Thus $[\sqrt{n!}]=k$.

\subsubsection{Remark}
\label{brocard_ramanujan_1st_remark}
All the known solutions follow the theorem as shown below.

(i) For $n=4$, $\sqrt{n!}\approx 4.898979486$. Then
$n!+1=[\sqrt{n!}]([\sqrt{n!}]+2)+1=(4\times 6)+1=25=5^2=m^2$.

(ii) For $n=5$, $\sqrt{n!}\approx 10.95445115$. Then
$n!+1=[\sqrt{n!}]([\sqrt{n!}]+2)+1=(10\times 12)+1=121=11^2=m^2$.

(iii) For $n=7$, $\sqrt{n!}\approx 70.992957397$. Then
$n!+1=[\sqrt{n!}]([\sqrt{n!}]+2)+1=(70\times 72)+1=5041=71^2=m^2$.

\subsubsection{Corollary}
\label{brocard_ramanujan_1st_corollary}
For the solution of equation (\ref{brocard_ramanujan_problem}),
$m-1=[\sqrt{n!}]$ should be an even natural number (i.e. $m$ should be
odd natural number) $k$ and $k(k+2)=n!$.

Example: (i) if $n=1$, then $k=[\sqrt{n!}]=1$ (odd number). Hence
no need to check solution of equation
(\ref{brocard_ramanujan_problem}).

(ii) if $n=2$, then $k=[\sqrt{n!}]=1$ (odd number). Hence no need to check solution of equation (\ref{brocard_ramanujan_problem}).

(iii) if $n=3$, then $k=[\sqrt{n!}]=2$, $k(k+2)=8>n!$. Hence equation (\ref{brocard_ramanujan_problem}) has no solution for $n=3$.

(iv) if $n=4$, then $k=[\sqrt{n!}]=4$, $k(k+2)=24=n!$, $n!+1=5^2$. Hence equation (\ref{brocard_ramanujan_problem}) has solution for $n=4$.

(v) if $n=5$, then $k=[\sqrt{n!}]=10$, $k(k+2)=120=n!$, $n!+1=11^2$. Hence equation (\ref{brocard_ramanujan_problem}) has solution for $n=5$.

(vi) if $n=6$, then $k=[\sqrt{n!}]=26$, $k(k+2)=728>n!$. Hence equation (\ref{brocard_ramanujan_problem}) has no solution for $n=6$.

(vii) if $n=7$, then $k=[\sqrt{n!}]=70$, $k(k+2)=5040=n!$, $n!+1=71^2$. Hence equation (\ref{brocard_ramanujan_problem}) has solution for $n=7$.

(viii) if $n=8$, then $k=[\sqrt{n!}]=26$, $k(k+2)=728>n!$. Hence equation (\ref{brocard_ramanujan_problem}) has no solution for $n=6$.

(ix) if $n=9$, then $k=[\sqrt{n!}]=602$, $k(k+2)=363608>n!$. Hence equation (\ref{brocard_ramanujan_problem}) has no solution for $n=9$.

(x) if $n=10$, then $k=[\sqrt{n!}]=1904$, $k(k+2)=3629024>n!$. Hence equation (\ref{brocard_ramanujan_problem}) has no solution for $n=10$.

(xi) if $n=11$, then $k=[\sqrt{n!}]=6317$ (odd number). Hence no need to
check solution of equation (\ref{brocard_ramanujan_problem}).

From the trend of the above example, we can predict the following
result: if $n$ is not a solution of Brocard-Ramanujan problem, then
$k(k+2)>n!$ where $k=[\sqrt{n!}]$. It is eventually proved in Theorem
\ref{brocard_ramanujan_4th_theorem}.

\subsubsection{Lemma}
\label{brocard_ramanujan_1st_lemma} 
For a given natural number $x$, let $x\neq z^2,~\text{for}~z\in
\mathbb{N}$ and $\sqrt{x}=y+\epsilon$, with
$y=[\sqrt{x}],~\epsilon=\sqrt{x}-[\sqrt{x}]~\text{i.e.}~0<\epsilon
<1$. Then $\epsilon (2y+\epsilon)$ is a natural number.

{\bf Proof}: Given that $\sqrt{x}=y+\epsilon$. Then
$x=y^2+2y\epsilon+\epsilon^2$. Or,
$x-y^2=2y\epsilon+\epsilon^2$. Since the left side is a natural
number, $\epsilon (2y+\epsilon)$ should also be a natural number.

\subsubsection{Remark}
\label{brocard_ramanujan_2nd_remark}
We can extend this idea for the $n$-th ($n\in \mathbb{N}$) root of
$x$. That is if $x\in \mathbb{N}$ with $x\neq z^n$ for $z,n\in
\mathbb{N}$ and $x^{\frac{1}{n}}=y+\epsilon$ having
$y=[x^{\frac{1}{n}}],~\epsilon=x^{\frac{1}{n}}-[x^{\frac{1}{n}}]~\text{i.e.}~0<\epsilon
<1$, then
$x-y^n=\binom{n}{1}y^{n-1}\epsilon+\binom{n}{1}y^{n-2}\epsilon^2+\cdots+\binom{n}{n}\epsilon^n$. Thus
$\binom{n}{1}y^{n-1}\epsilon+\binom{n}{2}y^{n-2}\epsilon^2+\cdots+\binom{n}{n}\epsilon^n$
is a natural number. For example, if $n=3$, then
$x^{\frac{1}{3}}=y+\epsilon$ with $y=[x^{\frac{1}{3}}],~0<\epsilon <1$
and $3y^{2}\epsilon+3y\epsilon^2+\epsilon^3$ is a natural number.

\subsection{Theorem}
\label{brocard_ramanujan_2nd_theorem}
In the Brocard-Ramanujan problem, let $\sqrt{n!}=k+\epsilon$ for the
natural number $n$ where $\sqrt{n!}\neq m~(\in
\mathbb{N}),~k=[\sqrt{n!}],~\epsilon=\sqrt{n!}-[\sqrt{n!}]$
i.e. $0<\epsilon <1$. Then $\epsilon (2k+\epsilon)$ must be a natural
number.

{\bf Proof}: Since $\sqrt{n!}\neq m~(\in \mathbb{N})$,
by the Lemma \ref{brocard_ramanujan_1st_lemma},
$n!=k^2+2k\epsilon+\epsilon^2$. Or,
$n!-k^2=2k\epsilon+\epsilon^2$. Thus $\epsilon (2k+\epsilon)$ is a
natural number. 

For example, (i) if $n=3$, then $\sqrt{n!}\approx 2.449489743$,
$k=2,~\epsilon\approx 0.449489743,~\epsilon (2k+\epsilon)=2$.

(ii) if $n=4$, then $\sqrt{n!}\approx 4.898979486$,
$k=4,~\epsilon\approx 0.898979486,~\epsilon (2k+\epsilon)=8$.

\subsection{Theorem}
\label{brocard_ramanujan_3rd_theorem}
Let for the equation (\ref{brocard_ramanujan_problem}),
$m-1=k=[\sqrt{n!}]$, $\sqrt{n!}=k+\epsilon$, where
$\epsilon=\sqrt{n!}-[\sqrt{n!}]$ i.e. $0<\epsilon <1$. Then for the
solution of (\ref{brocard_ramanujan_problem}), the natural number
$\epsilon (2k+\epsilon)$ must satisfy $\epsilon (2k+\epsilon)=2k$ and
conversely. Moreover, $k=\frac{\epsilon^2}{2(1-\epsilon)}$ and
$\epsilon$ should be distinct and monotonic increasing for all the
solutions of (\ref{brocard_ramanujan_problem}).

{\bf Proof}: Since $k=[\sqrt{n!}]$, we know from the Theorem
\ref{brocard_ramanujan_2nd_theorem} that $\epsilon (2k+\epsilon)$ is a
natural number. If $n$ is a solution of the equation
(\ref{brocard_ramanujan_problem}), then $k(k+2)=n!$ by Corollary
\ref{brocard_ramanujan_1st_corollary}. Or, $k(k+2)=(k+\epsilon)^2$.
\begin{eqnarray}
\text{Or,}~~~~\epsilon (2k+\epsilon)=2k.~~~~~~~~~~~~~~~~~~~~~~~~~~~~~~~~~~~~~~~~~~~~~~~~~~~~~~~~~~~~~~~~
\label{brocard_ramanujan_k_epsilon_expression}  
\end{eqnarray}
Hence the condition is necessary.

Conversely, let $\epsilon (2k+\epsilon)=2k$. Then
$k^2+2k\epsilon+\epsilon^2=k^2+2k$.

Or, $n!=k^2+2k$ (as $\sqrt{n!}=k+\epsilon$) i.e. $n!+1=(k+1)^2$. Thus
Brown's number $([\sqrt{n!}]+1,n)$ will be a solution of the equation
(\ref{brocard_ramanujan_problem}). Hence the condition is sufficient.
\begin{eqnarray}
\text{From (\ref{brocard_ramanujan_k_epsilon_expression}), we
  get}~~~~k=\frac{\epsilon^2}{2(1-\epsilon)}~~~~~~~~~~~~~~~~~~~~~~~~~~~~~~~~~~~~~~~~~~~~~~~~~~~~~~~~~~~~
\label{brocard_ramanujan_k_function_epsilon}
\end{eqnarray}
~~~~~~~~~~~and $\epsilon^2+2k\epsilon-2k=0$.

Or, $\epsilon=\frac{-2k\pm \sqrt{4k^2+8k}}{2}$.
\begin{equation}
\text{Then,}~~~~\epsilon=\sqrt{k^2+2k}-k.~~~~~~~~~~~~~~~~~~~~~~~~~~~~~~~~~~~~~~~~~~~~~~~~~~~~~~~~~~~~~~~~~~~~
\label{brocard_ramanujan_epsilon_function_k}
\end{equation}
It is easy to observe that $\epsilon>0$ and
$\epsilon<\sqrt{k^2+2k+1}-k=1$.

Let $\epsilon_k=f(k)=\sqrt{k^2+2k}-k$. 
\begin{equation}
\text{Then,}~~~~f^{\prime}(k)=\frac{2k+2}{2\sqrt{k^2+2k}}-1>\frac{k+1}{2\sqrt{k^2+2k+1}}-1=0~~~~~~~~~~~~~~~~~~~~~~~~~
\label{brocard_ramanujan_epsilon_mon_incr}
\end{equation}

Thus $\{f(k)\}=\{\epsilon_k\}$ is strictly monotonic increasing,
bounded above as $f(k)<1~\forall~k$ and also $\displaystyle
\lim_{k\rightarrow \infty}f(k)=\lim_{k\rightarrow
  \infty}\frac{2k}{\sqrt{k^2+2k}+k}=1$.

Hence, $\epsilon$ should be distinct and monotonic increasing for all
the solutions of equation (\ref{brocard_ramanujan_problem}).

Let us apply the theorem for the following first few cases.

(i) For $n=1$, $k=1,~\epsilon=0,~\frac{\epsilon^2}{2(1-\epsilon)}=0$
but $k(k+2)\neq n!$. Hence equation (\ref{brocard_ramanujan_problem})
has no solution for $n=1$.

(ii) For $n=2$, $k=1,~\epsilon\approx
0.414213562,~\frac{\epsilon^2}{2(1-\epsilon)}\approx
0.146446609~(\notin \mathbb{N}) <k$. Hence equation
(\ref{brocard_ramanujan_problem}) has no solution for $n=2$.

(iii) For $n=3$, $k=2,~\epsilon\approx
0.449489743,~\frac{\epsilon^2}{2(1-\epsilon)}\approx 0.183503419~(\notin
\mathbb{N}) <k$. Hence equation (\ref{brocard_ramanujan_problem}) has no
solution for $n=3$.

(iv) For $n=4$, $k=4,~\epsilon\approx
0.898979486,~\frac{\epsilon^2}{2(1-\epsilon)}=4$ and
$k(k+2)=24=n!$. Hence equation (\ref{brocard_ramanujan_problem}) has
solution for $n=4$.

(v) For $n=5$, $k=10,~\epsilon\approx
0.95445115,~\frac{\epsilon^2}{2(1-\epsilon)}=10$ and
$k(k+2)=120=n!$. Hence equation (\ref{brocard_ramanujan_problem}) has
solution for $n=5$. 

(vi) For $n=6$, $k=26,~\epsilon\approx
0.83281573,~\frac{\epsilon^2}{2(1-\epsilon)}\approx 2.07430412~(\notin
\mathbb{N}) <k$ and also $k(k+2)\neq n!$. Hence equation
(\ref{brocard_ramanujan_problem}) has no solution for $n=6$.

(vii) For $n=7$, $k=70,~\epsilon\approx
0.9929573972,~\frac{\epsilon^2}{2(1-\epsilon)}=70$ and
$k(k+2)=5040=n!$. Hence equation (\ref{brocard_ramanujan_problem}) has
solution for $n=7$. 

(viii) For $n=8$, $k=200,~\epsilon\approx
0.7984064,~\frac{\epsilon^2}{2(1-\epsilon)}\approx 1.581033892~(\notin
\mathbb{N}) <k$ and also $k(k+2)\neq n!$. Hence equation
(\ref{brocard_ramanujan_problem}) has no solution for $n=8$.

(ix) For $n=9$, $k=602,~\epsilon\approx
0.3952191,~\frac{\epsilon^2}{2(1-\epsilon)}\approx 0.1291361398~(\notin
\mathbb{N})<k$ and also $k(k+2)\neq n!$. Hence equation
(\ref{brocard_ramanujan_problem}) has no solution for $n=9$.

(x) For $n=10$, $k=1904,~\epsilon\approx
0.940944,~\frac{\epsilon^2}{2(1-\epsilon)}\approx 7.496063447~(\notin
\mathbb{N})<k$ and also $k(k+2)\neq n!$. Hence equation
(\ref{brocard_ramanujan_problem}) has no solution for $n=10$.

(xi) For $n=11$, $k=6371,~\epsilon\approx
0.974359,~\frac{\epsilon^2}{2(1-\epsilon)}\approx 18.51278095~(\notin
\mathbb{N})<k$ and also $k(k+2)\neq n!$. Hence equation
(\ref{brocard_ramanujan_problem}) has no solution for $n=11$.

\subsubsection{Remark}
\label{brocard_ramanujan_3rd_remark}
(i) The function $f(k)$ has no optimum value in $0<\epsilon<1$. For the
known solutions of Brocard-Ramanujan problem (c.f. equation
(\ref{brocard_ramanujan_problem})), (i) $\epsilon\approx 0.898979486$
for $n=3$ (ii) $\epsilon\approx 0.95445115$ for $n=5$ (iii)
$\epsilon\approx 0.9929573972$ for $n=7$. The next value of $\epsilon$
should be more than $\epsilon\approx 0.9929573972$ as $\epsilon$
should be distinct and monotonic increasing for all the solutions of
equation (\ref{brocard_ramanujan_problem}). However, $\epsilon$ is
neither monotonic increasing nor monotonic decreasing for all natural
numbers $n$ (c.f. examples of Theorem
\ref{brocard_ramanujan_3rd_theorem}). Also, $\epsilon$ is a function
of $k,~k=[\sqrt{n!}]$ in the form of the equation
(\ref{brocard_ramanujan_epsilon_function_k}) and according to Bruce
Berndt and William Galway \cite{Berndt}, there is no new solution of
the problem based on computations up to $n=10^9$. If we consider a new
solution of Brocard-Ramanujan problem as $n\ge 10^5$, then value of
$\epsilon$ should be more than $\epsilon\approx 0.999\cdots905915$
({\bf digit 0 is coming after 228287 numbers of 9 digit, which takes
  more than 66 pages in LibreOffice Writer}) for $k=[\sqrt{10^5!}]$
according to the calculation in MATHEMATICA 8. If we consider $n\geq
10^9$ as new solution of Brocard-Ramanujan problem, then it is obvious
that the numbers of 9 digit in the form $\epsilon=0.999\cdots$ will be
much higher than the corresponding case of $n\geq 10^5$, but I am
unable to calculate the numbers of 9 digit in the value of $\epsilon$
in my personal laptop (with 8GB Ram) using MATHEMATICA 8. Hence, it suggests that there will be no new solution of the problem.

(ii) For the non-solution of (\ref{brocard_ramanujan_problem}), what will be the expression of the natural number
$\epsilon (2k+\epsilon)?$. The said expression is provided in theorem \ref{brocard_ramanujan_13th_theorem}.

\subsubsection{Corollary}
\label{brocard_ramanujan_2nd_corollary}
For the solution of Brocard-Ramanujan problem, we have
$k=\frac{\epsilon^2}{2(1-\epsilon)}$. Hence,
$\{g(\epsilon)\}=\{k_\epsilon\}$ is also monotonic increasing.

\subsection{Theorem}
\label{brocard_ramanujan_4th_theorem}
For any natural number $n$, $n!\le
k(k+2),~\text{where}~k=[\sqrt{n!}]$ and if $n$ is not a solution of
Brocard-Ramanujan problem (c.f. equation
(\ref{brocard_ramanujan_problem})), then $n!<k(k+2)$.

{\bf Proof}: For any natural number $n$, $\sqrt{n!}=k+\epsilon$ with
$k=[\sqrt{n!}]$ and $0\le \epsilon<1$. Then
$n!=k^2+2k\epsilon+\epsilon^2<k^2+2k+1$ as $\epsilon<1$. We know that
$n!,~k^2$ and $2k\epsilon+\epsilon^2$ (by Theorem
\ref{brocard_ramanujan_3rd_theorem}) are integers. Hence $n!\le
k(k+2),~\forall~n\in \mathbb{N}$.

We also know from Corollary \ref{brocard_ramanujan_1st_corollary} that
for the solution of Brocard-Ramanujan problem, we have
$n!=k(k+2)$. Thus if $n$ is not a soluton of Brocard-Ramanujan problem
(c.f. equation (\ref{brocard_ramanujan_problem})), then $n!<k(k+2)$.

\subsubsection{Corollary}
\label{brocard_ramanujan_3rd_corollary}
If natural number $n$ is not a solution of Brocard-Ramanujan problem
(c.f. equation (\ref{brocard_ramanujan_problem})), then
$2k\epsilon+\epsilon^2<2k$.

\subsection{Theorem}
\label{brocard_ramanujan_5th_theorem}
The Brocard-Ramanujan problem (c.f. equation
(\ref{brocard_ramanujan_problem})) has only finitely many solutions.

{\bf Proof}: From Theorem \ref{brocard_ramanujan_3rd_theorem}, we know
that for the solution of (\ref{brocard_ramanujan_problem}), the
natural number $2k\epsilon+\epsilon^2$ must satisfy
$2k\epsilon+\epsilon^2=2k$. Our endeavor should be to find all possible $k$
and corresponding $\epsilon$ such that $2k\epsilon+\epsilon^2=2k$ for
the Brocard-Ramanujan problem. Moreover, if $n$ is not a solution of
Brocard-Ramanujan problem, then $2k\epsilon+\epsilon^2<2k$ by
Corollary \ref{brocard_ramanujan_3rd_corollary}.

Let $f(\epsilon,k)=2k\epsilon+\epsilon^2$. Then by Taylor's expansion
of two variables with respect to $(\epsilon_0,k_0)$, we obtain
$f(\epsilon,k)=f(\epsilon_0,k_0)+(\epsilon-\epsilon_0)\frac{\partial
	f}{\partial \epsilon}(\epsilon_1,k_1)+(k-k_0)\frac{\partial
	f}{\partial
	k}(\epsilon_1,k_1)=2k_0\epsilon_0+\epsilon_0^2+(\epsilon-\epsilon_0)(2\epsilon_1+2k_1)+(k-k_0)(2\epsilon_1)$
with $k_0=[\sqrt{7!}]=70,~\epsilon_0=\sqrt{7!}-[\sqrt{7!}]\approx
0.9929573972$ and for some $\epsilon_1$ and $k_1$ satisfying
$\epsilon_0<\epsilon_1<\epsilon<1,~k_0<k_1<k$.

Thus
$f(\epsilon,k)=140+(\epsilon-\epsilon_0)(2\epsilon_1+2k_1)+2k+2k(\epsilon_1-1)-140\epsilon_1=2k+140(1-\epsilon_1)+(\epsilon-\epsilon_0)(2\epsilon_1+2k_1)-2k(1-\epsilon_1)<2k$
if
$140(1-\epsilon_1)+(\epsilon-\epsilon_0)(2\epsilon_1+2k_1)-2k(1-\epsilon_1)<0$
i.e. $k>\frac{140(1-\epsilon_1)+2(\epsilon-\epsilon_0)(\epsilon_1+k_1)}{2(1-\epsilon_1)}=70+\frac{(\epsilon-\epsilon_0)(\epsilon_1+k_1)}{(1-\epsilon_1)}$.

Thus $f(\epsilon,k)<2k$ for all
$k>70+\frac{(1-\epsilon_0)(\epsilon_1+k_1)}{(1-\epsilon_1)}$.

Then, the Brocard-Ramanujan problem (c.f. equation
(\ref{brocard_ramanujan_problem})) has no solutions if
$k>70+\frac{(1-\epsilon_0)(\epsilon_1+k_1)}{(1-\epsilon_1)}$.

Hence, the Brocard-Ramanujan problem has only finitely many solutions.

\subsection{Theorem}
\label{brocard_ramanujan_6th_theorem}
Let for the equation (\ref{brocard_ramanujan_problem}),
$m-1=k=[\sqrt{n!}]$, $\sqrt{n!}=k+\epsilon$, where
$\epsilon=\sqrt{n!}-[\sqrt{n!}]$ i.e. $0<\epsilon <1$. Then for the
solution of (\ref{brocard_ramanujan_problem}),
$n!=\frac{\epsilon^2(2-\epsilon)^2}{4(1-\epsilon)^2}$ i.e. the
solutions $(m, n!)$ will be in the form
$(\frac{\epsilon^2}{2(1-\epsilon)}+1,\frac{\epsilon^2(2-\epsilon)^2}{4(1-\epsilon)^2})$. Moreover,
for all the solutions of (\ref{brocard_ramanujan_problem}), $x=n!$ and
$y=m-1=k=[\sqrt{n!}]$ should satisfy
\begin{eqnarray}
y^4+4y^3+2xy^2+4y^2+4xy-3x^2=0,~y^3+3y^2+2y-xy-x=0.
\end{eqnarray}

{\bf Proof}: For the solution of Brocard-Ramanujan problem, we have
$k=\frac{\epsilon^2}{2(1-\epsilon)}$. Then
$\sqrt{n!}=k+\epsilon=\frac{\epsilon(2-\epsilon)}{2(1-\epsilon)}$.
\begin{eqnarray}
\text{Thus,}~~n!=\frac{\epsilon^2(2-\epsilon)^2}{4(1-\epsilon)^2}.~~~~~~~~~~~~~~~~~~~~~~~~~~~~~~~~~~~~~~~~~~~~~~~~~~~~~~~~~
\label{brocard_ramanujan_n!_epsilon}
\end{eqnarray}
Or, $4n!(1-2\epsilon+\epsilon^2)=\epsilon^2(4-4\epsilon+\epsilon^2)$.

Or,
$4n!\{1-2(\sqrt{n!}-[\sqrt{n!}])+(\sqrt{n!}-[\sqrt{n!}])^2\}=\{(\sqrt{n!}-[\sqrt{n!}])^2(4-4(\sqrt{n!}-[\sqrt{n!}])+(\sqrt{n!}-[\sqrt{n!}])^2\}$,
where $\epsilon=\sqrt{n!}-[\sqrt{n!}]$.

Or, $4n!\{(1+2[\sqrt{n!}]+n!+[\sqrt{n!}]^2\}-8n!\sqrt{n!}\{1+[\sqrt{n!}]\}=\{(4n!+4[\sqrt{n!}]^2+12n![\sqrt{n!}]+n!^2+6n![\sqrt{n!}]^2+4[\sqrt{n!}]^3+[\sqrt{n!}]^4)-4n!\sqrt{n!}(2[\sqrt{n!}]+n!+3[\sqrt{n!}]^2+n![\sqrt{n!}]+[\sqrt{n!}]^3)\}$.

Comparing both sides, we get $4n!\{(1+2[\sqrt{n!}]+n!+[\sqrt{n!}]^2\}=\{(4n!+4[\sqrt{n!}]^2+12n![\sqrt{n!}]+n!^2+6n![\sqrt{n!}]^2+4[\sqrt{n!}]^3+[\sqrt{n!}]^4),~2n!+2n![\sqrt{n!}]=2[\sqrt{n!}]+n!+3[\sqrt{n!}]^2+n![\sqrt{n!}]+[\sqrt{n!}]^3$.

\begin{eqnarray}
\text{Thus,}~~y^4+4y^3+2xy^2+4y^2+4xy-3x^2=0,
\label{brocard_ramanujan_n_k_relation_1}\\
y^3+3y^2+2y-xy-x=0;
\label{brocard_ramanujan_n_k_relation_2}
\end{eqnarray}
where $x=n!,~y=[\sqrt{n!}]$.

All the known solutions of Brocard-Ramanujan problem satisfy the
equations (\ref{brocard_ramanujan_n_k_relation_1}) and
(\ref{brocard_ramanujan_n_k_relation_2}) as shown below.

(i) For $n=4$, $x=24,~y=4$. Then
$y^4+4y^3+2xy^2+4y^2+4xy-3x^2=0$ and $y^3+3y^2+2y-xy-x=0$.

(ii) For $n=5$, $x=120,~y=10$. Then $y^4+4y^3+2xy^2+4y^2+4xy-3x^2=0$ and $y^3+3y^2+2y-xy-x=0$.

(iii) For $n=7$, $x=5040,~y=70$. Then $y^4+4y^3+2xy^2+4y^2+4xy-3x^2=0$ and $y^3+3y^2+2y-xy-x=0$.

\subsubsection{Remark}
\label{brocard_ramanujan_4th_remark}
It is shown that all solutions of Brocard-Ramanujan problem satisfy the
equations (\ref{brocard_ramanujan_n_k_relation_1}) and
(\ref{brocard_ramanujan_n_k_relation_2}). By B$\acute{e}$zout theorem,
(c.f. Sturmfels \cite{ Sturmfels}) for a general square polynomial
system
$f_1(x_1,x_2,\cdots,x_n)=0,~f_2(x_1,x_2,\cdots,x_n)=0,~\cdots,~f_n(x_1,x_2,\cdots,x_n)=0$
having $d_1,d_2,\cdots,d_n$ as the degrees of $f_1,f_2,\cdots,f_n$
respectively; the number of its isolated zeros in $\mathbb{C}^n$,
counting multiplicities, does not exceed the number $d=d_1d_2\cdots
d_n$ except the system has an infinite number of zeros.

However, in-built functions `Solve', `Solve with integer', `NSolve',
`NSolve with real', `NSolve with integer', `FindRoot' are unable to
provide any solutions of Brocard-Ramanujan problem. To solve by hands,
the following approaches are considered.

{\bf Approach 1}: If we solve equation
(\ref{brocard_ramanujan_n_k_relation_2}), then we get
$y=-1,~x=y(y+2)$. For $y=-1$, we get from equation
(\ref{brocard_ramanujan_n_k_relation_1}) as $3x^2+2x-1=0$
i.e. $x=-1,\frac{1}{3}$.

For $x=y(y+2)$, we get from equation
(\ref{brocard_ramanujan_n_k_relation_1}) as
$(4y+2y^2)y(y+2)-3y^2(y+2)^2+4y^2+4y^3+y^4=0$. Or,
$2y^2(y+2)^2-3y^2(y+2)^2+y^2(y+2)^2=0$ Or, $0=0$, which gives no
explicit solution in this case. However, we can get three known
solutions of the Brocard-Ramanujan problem as
$(24,4),~(120,10),~(5040,70)$ and infinite number of integer solutions
in the form $(c,d)$ with $d<0,~d\in\mathbb{Z}$, $c=d^2-2d$. First ten
of these solutions are $(-1,-1),~(0,-2),~(3,-3),~(8,-4),~(15,-5),$
$(24,-6),~(35,-7),~(48,-8),~(63,-9),~(80, -10)$.

{\bf Approach 2}: From equation
(\ref{brocard_ramanujan_n_k_relation_1}), we have
$3x^2-2x(y^2+2y)-(4y^2+4y^3+y^4)=0$. Or, $x=\frac{2y(y+2)\pm
  \sqrt{4y^2(y+2)^2+12(4y^2+4y^3+y^4)}}{6}=\frac{y(y+2)\pm
  2y(y+2)}{3}=y(y+2),~-\frac{y(y+2)}{3}$. 

If $x=y(y+2)$, then the equation
(\ref{brocard_ramanujan_n_k_relation_2}) gives
$-(y+1)y(y+2)+y^3+3y^2+2y=0$. Or, $-y(y+1)y(y+2)+y(y+1)y(y+2)=0$. Or,
$0=0$, which gives no explicit solution in this case. But, we can find
solutions as $(24,4),~(120,10),~(5040,70)$ and $(c,d)$ with
$d<0,~d\in\mathbb{Z}$, $c=d^2-2d$.

If $x=-\frac{y(y+2)}{3}$, then the equation
(\ref{brocard_ramanujan_n_k_relation_2}) gives
$-\frac{(y+1)y(y+2)}{3}+y^3+3y^2+2y=0$. Or,
$-\frac{(y+1)y(y+2)}{3}+(y+1)y(y+2)=0$. Or,
$\frac{2(y+1)y(y+2)}{3}=0$. Or, $y=0,-1,-2$. Then $x=0,\frac{1}{3},0$
respectively.

{\bf Approach 3}: From equation
(\ref{brocard_ramanujan_n_k_relation_1}), we have
$y^4+4y^3+(2x+4)y^2+4xy-3x^2=0$. 

\begin{eqnarray}
\text{Or,}~~(y^2+2y+\lambda)^2-(my+n)^2=0
\label{brocard_ramanujan_ferrari_method}
\end{eqnarray}
with $4+2\lambda -m^2=2x+4,~4\lambda-2mn=4x,~\lambda^2-n^2=-3x^2$
i.e. $m^2=2(\lambda-x),~mn=2(\lambda-x),~n^2=\lambda^2+3x^2$. Then $m^2=mn$ i.e. $m=0,n$.

For $m=0$, $\lambda=x$, $n^2=4x^2$ i.e. $n=\pm 2x$. Then from equation
(\ref{brocard_ramanujan_ferrari_method}), we have $(y^2+2y+x)^2-(\pm
2x)^2=0$. Or, $(y^2+2y+3x)(y^2+2y-x)=0$. Or,
$x=y(y+2),~-\frac{y(y+2)}{3}$. If $x=y(y+2)$, then the equation
(\ref{brocard_ramanujan_n_k_relation_2}) gives
$-(y+1)y(y+2)+y^3+3y^2+2y=0$. Or, $-y(y+1)y(y+2)+y(y+1)y(y+2)=0$. Or,
$0=0$, which gives no explicit solution in this case and we can find
solutions as $(24,4),~(120,10),~(5040,70)$ and $(c,d)$ with
$d<0,~d\in\mathbb{Z}$, $c=d^2-2d$.

If $x=-\frac{y(y+2)}{3}$, then the equation
(\ref{brocard_ramanujan_n_k_relation_2}) gives
$-\frac{(y+1)y(y+2)}{3}+y^3+3y^2+2y=0$. Or,
$-\frac{(y+1)y(y+2)}{3}+y(y+1)y(y+2)=0$. Or,
$\frac{2(y+1)y(y+2)}{3}=0$. Or, $y=0,-1,-2$.

For $m=n$, $n^2=2(\lambda-x)$ and $n^2=\lambda^2+3x^2$
i.e. $2(\lambda-x)=\lambda^2+3x^2$. Or,
$\lambda^2-2\lambda+2x+3x^2=0$. Or, $\lambda=\frac{2\pm
  \sqrt{4-4(3x^2+2x)}}{2}=\frac{1\pm \sqrt{1-(3x^2+2x)}}{2}$. We can
see that $1-(3x^2+2x)<0$ for $x>1$.

Hence there are solutions
$(0,0),~(\frac{1}{3},-1),~(24,4),~(120,10),~(5040,70)$ and $(c,d)$ with
$d<0,~d\in\mathbb{Z}$, $c=d^2-2d$ of equations
(\ref{brocard_ramanujan_n_k_relation_1}) and
(\ref{brocard_ramanujan_n_k_relation_2}); among them only three
solutions $(24,4),~(120,10),$ $(5040,70)$ are of the Brocard-Ramanujan
problem which are already known.

\subsection{Theorem}
\label{brocard_ramanujan_7th_theorem}
For the fourth solution of the Brocard-Ramanujan
problem (\ref{brocard_ramanujan_problem}), $m=\alpha \times 10^{\beta}\pm 1$ for some $\alpha,~\beta\in \mathbb{N}$ (i) with $\beta\geq 24999$ when $n\geq 10^5$ and (ii) with $\beta\geq 249999998$ when $n\geq 10^9$. 

{\bf Proof}: We know that $n!=(m-1)(m+1)$ with gcd $d=(m-1,m+1)=2$ and $2^\beta_1/n!$ but $2^{\beta_1+1}$ not divide $n!$ where $\beta_1=[\frac{n}{2}]+[\frac{n}{4}]+[\frac{n}{8}]+\cdots$. Then $2^{\beta_1-1}/(m-1)$ or $2^{\beta_1-1}/(m+1)$ but not both. Again, $5^{\beta_3}/n!$ but $5^{\beta_3+1}$ not divide $n!$ where $\beta_3=[\frac{n}{5}]+[\frac{n}{25}]+[\frac{n}{125}]+\cdots$. and $5^{\beta_3}/(m-1)$ or $5^{\beta_3}/(m+1)$ but not both. 

If $2/(m-1),~5^{\beta_3}/(m-1)$ and $2^{\beta_1-1}/(m+1)$, then $10/(m-1)$ but $100$ not divide $(m-1)$ and $10$ does not divide $m+1$ which contradicts $100/n!$ for $n\geq 10$ as $d=2$. Also for $2/(m+1),~5^{\beta_3}/(m+1)$ and $2^{\beta_1-1}/(m-1)$, then $10/(m+1)$ but $100$ not divide $(m+1)$ and $10$ does not divide $m-1$ which contradicts $100/n!$ for $n\geq 10$ as $d=2$. Hence if $2^{\beta_1-1}/(m-1)$, them $5^{\beta_3}/(m-1)$ and if $2^{\beta_1-1}/(m+1)$, them $5^{\beta_3}/(m+1)$.

(i) For $n=10^5$, we can calculate as $\beta_3=24999$ and $\beta_1>\beta_3$. Thus $m$ will be in the form $m=\alpha \times 10^{\beta}\pm 1$ for some $\alpha,~\beta\in \mathbb{N}$ with $\beta\geq 24999$ when $n\geq 10^5$.

(ii) For $n=10^9$, we can evaluate as $\beta_3=249999998$ and $\beta_1>\beta_3$. Thus $m$ will be in the form $m=\alpha \times 10^{\beta}\pm 1$ for some $\alpha,~\beta\in \mathbb{N}$ with $\beta\geq 249999998$ when $n\geq 10^9$.

\subsection{Theorem}
\label{brocard_ramanujan_8th_theorem}
For the fourth solution of the Brocard-Ramanujan
problem (\ref{brocard_ramanujan_problem}), $m=2^{\beta_1-1}5^{\beta_3}m_3+1$ or $m=2^{\beta_1-1}5^{\beta_3}m_3-1$, $m_3$ is odd.

{\bf Proof}: For the Brocard-Ramanujan
problem (\ref{brocard_ramanujan_problem}), we have $n!=m^2-1=(m-1)(m+1),$ where $m$ is odd.

Case I: $m=4m_1+1$. Then $n!=4\times 2m_1\times(2m_1+1)$.

Case I(i): $m_1=2m_2$. Then $m=8m_2+1,~n!=4\times 4m_2\times(4m_2+1)=4\times 2^{\beta_1-2}5^{\beta_3}m_3\times (2^{\beta_1-2}5^{\beta_3}m_3+1)$ using theorem (\ref{brocard_ramanujan_7th_theorem}), where $m_2=2^{\beta_1-4}5^{\beta_3}m_3,~m=2^{\beta_1-1}5^{\beta_3}m_3+1$ and $m_3$ is odd.

Case I(ii): $m_1=2m_2+1$. Then $m=8m_2+5,~n!=4\times 2(2m_2+1)\times (4m_2+3)$. Then here $\beta_1=3$ and so no fourth solution possible in this case as $\beta_1\geq 999999987$ when $n\geq 10^9$. 

Case II: $m=4m_1+3$. Then $n!=4\times (2m_1+1)\times 2(m_1+1)$.

Case II(i): $m_1=2m_2$. Then $m=8m_2+5,~n!=4\times (4m_2+1)\times 2(2m_2+1)$. Then here $\beta_1=3$ and so no fourth solution possible in this case as $\beta_1\geq 999999987$ when $n\geq 10^9$.

Case II(ii): $m_1=2m_2+1$. Then $m=8m_2+7,~n!=4\times (4m_2+3)\times 4(m_2+1)=4\times (2^{\beta_1-2}5^{\beta_3}m_3-1)\times 2^{\beta_1-2}5^{\beta_3}m_3$ using theorem (\ref{brocard_ramanujan_7th_theorem}), where $m_2+1=2^{\beta_1-4}5^{\beta_3}m_3,~m=2^{\beta_1-1}5^{\beta_3}m_3-1$ and $m_3$ is odd.

\subsection{Theorem}
\label{brocard_ramanujan_9th_theorem}
For all the solutions of the Brocard-Ramanujan
problem (\ref{brocard_ramanujan_problem}), $m$ should be in the form of $m=5+6t$ and $m=11+30t$ for all the solutions of the Brocard-Ramanujan problem (\ref{brocard_ramanujan_problem}) except first, while $m=71+5040t$ for all the solutions starting from third.

{\bf Proof}: From the result \cite{Mapa}, we know that if $d_1,~d_2\in\mathbb{N}$ and $e_1,~e_2\in\mathbb{Z}$; then $y=e_1$ (mod $d_1$), $y=e_2$ (mod $d_2$) have a simultaneous solution iff gcd $(d_1,d_2)/(e_1-e_2)$ and if this condition be satisfied; the solution is unique module lcm $(d_1,d_2)$.

For the solution of the Brocard-Ramanujan
problem (\ref{brocard_ramanujan_problem}), we have $n!=m^2-1$. Then $m=1$ (mod $\frac{k}{2})$ and $m=-1$ (mod $\frac{k+2}{2})$ with gcd $(\frac{k}{2},\frac{k+2}{2})=1$.

Case I:  For the Brown's number $(5,4)$, $m=5,~k=4$, $m=1$ (mod $2$) and $m=-1$ (mod $3$) and the solution is $f_t=m+6t=5+6t$. When $t=1,~g_1=11$ and also $11=1$ (mod $2$), $11=-1$ (mod $3$) giving $m$ for 2nd solution of (\ref{brocard_ramanujan_problem}). For $t=1,~g_11=71$ and also $71=1$ (mod $2$), $71=-1$ (mod $3$) which provides $m$ for 3rd solution of (\ref{brocard_ramanujan_problem}).

Case II:  For the Brown's number $(11,5)$, $m=11,~k=10$, $m=1$ (mod $5$) and $m=-1$ (mod $6$) and the solution is $f_t=m+30t=11+30t$. When $t=2,~g_2=71$ and also $71=1$ (mod $5$), $71=-1$ (mod $6$) giving $m$ for 3rd solution of (\ref{brocard_ramanujan_problem}). 

Case III:  For the Brown's number $(71,7)$, $m=71,~k=70$, $m=1$ (mod $35$) and $m=-1$ (mod $36$) and the solution is $f_t=m+5040t=71+5040t$.

\subsubsection{Note}
In all three cases, $n$ can be evaluated as $n!=m^2-1$.

\subsection{Theorem}
\label{brocard_ramanujan_10th_theorem}
If $\sqrt{n!}=k+\epsilon$ for solution of Brocard-Ramanujan problem (c.f. equation
(\ref{brocard_ramanujan_problem})), then necessary condition is $k$ divides $n!$, while for sufficient condition $\frac{n!}{k}=k+2$. If $2k\epsilon+\epsilon^2=k$, then there will be no other solution of (\ref{brocard_ramanujan_problem}) in this case. But computation suggests that $n=3$ may be only one example for which $2k\epsilon+\epsilon^2=k$. If $2k\epsilon+\epsilon^2=k$, then $\epsilon<0.5$ and $\epsilon,~k$ are strictly monotonic increasing for the solutions of $n!=k(k+1)$. 

{\bf Proof}: $n!=k(k+1)$ If $n$ is a solution of the equation
(\ref{brocard_ramanujan_problem}), then $k(k+2)=n!$ by Corollary
\ref{brocard_ramanujan_1st_corollary}. Thus necessary condition of (\ref{brocard_ramanujan_problem}) is that $k$ divides $n!$. But it is not sufficient condition as we need equation $\frac{n!}{k}=k+2$. Then $n!+1=(k+1)^2$. Hence the sufficient condition is $\frac{n!}{k}=k+2$.

If $2k\epsilon+\epsilon^2=k$, then $n!=k^2+2k\epsilon+\epsilon^2=k(k+1)$ and $\frac{n!}{k}=k+1$. Then there will be no other solution of (\ref{brocard_ramanujan_problem}) in this case. 

If $2k\epsilon+\epsilon^2=k$, it is clear that $\epsilon<0.5$ and $\epsilon\longrightarrow 0.5$ as $n\longrightarrow \infty$. Also $\epsilon=\sqrt{k^2+k}-k$. Or, $\frac{d\epsilon}{dk}=\frac{2k+1}{2\sqrt{k^2+k}}-1$ and $\left(\frac{2k+1}{2\sqrt{k^2+k}}\right)^2=\frac{4k^2+4k+1}{4k^2+4k}>1$. Then $\frac{d\epsilon}{dk}>0$. Hence $\epsilon$ is strictly monotonic increasing for the solutions of $n!=k(k+1)$. Again, $k=\frac{\epsilon^2}{1-2\epsilon}$, $\frac{dk}{d\epsilon}=\frac{2\epsilon(1-\epsilon)}{(1-2\epsilon)^2}>0$. So, $k$ is strictly monotonic increasing for the solutions of $n!=k(k+1)$

We know that for $n=3$, $3!=6=2(2+1),~k=[\sqrt{3!}]=2,~\epsilon\approx 0.4494897428,~2k\epsilon+\epsilon^2=2=k$.

If $n=17$, $k=18859677$. when $\epsilon$ satisfies $2k\epsilon+\epsilon^2=k$, then $\epsilon\approx 0.499999993$. But actual value of $\epsilon$ for $n=17$ is $\epsilon\approx 0.306253$.

\subsection{Theorem}
\label{brocard_ramanujan_11th_theorem}
If $\sqrt{n!}=k+\epsilon$ for solution of Brocard-Ramanujan problem (c.f. equation
(\ref{brocard_ramanujan_problem})), then $(\epsilon,k)$ lies on hyperbola.

{\bf Proof}: By the Theorem \ref{brocard_ramanujan_3rd_theorem} for the
solution of (\ref{brocard_ramanujan_problem}), the natural number $\epsilon (2k+\epsilon)$ must satisfy $\epsilon (2k+\epsilon)=2k$ i.e. $\epsilon^2+2\epsilon k-2k=0$, which is a hyperbola with center $(1,-1)$. Let us consider the translation $\epsilon=\epsilon_1+1,~k=k_1-1$. Then the hyperbola is shifted to $\epsilon_1^2+2\epsilon_1k_1+1=0$. Using the rotation $\epsilon_1=\epsilon_2\cos \theta-k_2\sin \theta,~k_1=\epsilon_2\sin \theta+k_2\cos \theta$, we get the transformed hyperbola as $(\cos^2\theta+2\sin\theta\cos\theta)\epsilon_2^2+(\sin^2\theta-2\sin\theta\cos\theta)k_2^2+(-2\sin \theta\cos\theta+2\cos^2\theta-2\sin^2\theta)\epsilon_2k_2+1=0$. For the value of $\theta$, we consider coefficient of $\epsilon_2k_2$ as zero. Then $-2\sin \theta\cos\theta+2\cos^2\theta-2\sin^2\theta=0$. Or, $\tan^2\theta+\tan\theta-1=0$. Or, $\tan\theta=\frac{-1\pm\sqrt{5}}{2}$. 

For $\tan\theta=\frac{\sqrt{5}-1}{2},~\sin\theta=\frac{\sqrt{5}-1}{\sqrt{10-2\sqrt{5}}},~\cos\theta=\frac{2}{\sqrt{10-2\sqrt{5}}}$. Then the transformed hyperbola is $\frac{3-\sqrt{5}}{\sqrt{5}-1}k_2^2-\frac{2}{\sqrt{5}-1}\epsilon_2^2=1,$ where $\epsilon=1+\frac{2}{\sqrt{10-2\sqrt{5}}}\epsilon_2-\frac{\sqrt{5}-1}{\sqrt{10-2\sqrt{5}}}k_2$ and $k=-1+\frac{\sqrt{5}-1}{\sqrt{10-2\sqrt{5}}}\epsilon_2+\frac{2}{\sqrt{10-2\sqrt{5}}}k_2$. 

For $\tan\theta=\frac{-\sqrt{5}-1}{2},~\sin\theta=\frac{\sqrt{5}+1}{\sqrt{10+2\sqrt{5}}},~\cos\theta=-\frac{2}{\sqrt{10+2\sqrt{5}}}$ or $\sin\theta=-\frac{\sqrt{5}+1}{\sqrt{10+2\sqrt{5}}},~\cos\theta=\frac{2}{\sqrt{10+2\sqrt{5}}}$. Then the transformed hyperbola is $\frac{2}{\sqrt{5}+1}\epsilon_2^2-\frac{3+\sqrt{5}}{\sqrt{5}+1}k_2^2=1,$ where $\epsilon=1-\frac{2}{\sqrt{10+2\sqrt{5}}}\epsilon_2-\frac{\sqrt{5}+1}{\sqrt{10+2\sqrt{5}}}k_2$ and $k=-1+\frac{\sqrt{5}+1}{\sqrt{10+2\sqrt{5}}}\epsilon_2-\frac{2}{\sqrt{10+2\sqrt{5}}}k_2$ or $\epsilon=1+\frac{2}{\sqrt{10+2\sqrt{5}}}\epsilon_2+\frac{\sqrt{5}+1}{\sqrt{10+2\sqrt{5}}}k_2$ and $k=-1-\frac{\sqrt{5}+1}{\sqrt{10+2\sqrt{5}}}\epsilon_2+\frac{2}{\sqrt{10+2\sqrt{5}}}k_2$ respectively. 

\subsubsection{Note}
(i) $2k=\epsilon^2+2k\epsilon<3k\epsilon$ i.e. $\epsilon>\frac{2}{3}$ for the solution of (\ref{brocard_ramanujan_problem}).

(ii) $2k=\epsilon^2+2k\epsilon<\frac{k\epsilon}{2}+2k\epsilon$ as $k\geq 2$ for $n\geq 3$ i.e. $\epsilon>\frac{4}{5}$ for the solution of (\ref{brocard_ramanujan_problem}).

(iii) $2k=\epsilon^2+2k\epsilon<1+2k\epsilon$ i.e. $\epsilon>\frac{2k-1}{2k}$. For $n=4$, $k=4$ and hence $\epsilon>\frac{7}{8}$ for the solution of (\ref{brocard_ramanujan_problem}) if $n\geq 4$.

(iv) $2k=\epsilon^2+2k\epsilon<\epsilon+2k\epsilon$ i.e. $\epsilon>\frac{2k}{2k+1}$. For $n=4$, $k=4$ and hence $\epsilon>\frac{8}{9}$ for the solution of (\ref{brocard_ramanujan_problem}) if $n\geq 4$. It can be noted that $\frac{2k}{2k+1}$ better approximation of $\epsilon$ than $\frac{2k-1}{2k}$.

(v) $\epsilon^2+2k\epsilon\geq (1+2k)\left(\epsilon^{2k}\epsilon^2\right)^{\frac{1}{2k+1}}$. Or, $(2k+1)\epsilon^{\frac{2k+2}{2k+1}}\leq \epsilon^2+2k\epsilon<g(2k+1)$ for some $0<g<1$. Then $\epsilon<g^{\frac{2k+1}{2k+2}}$. Also $2k=\epsilon^2+2k\epsilon<g(2k+1)$ i.e. $g>\frac{2k}{2k+1}$.

\subsection{Theorem}
\label{brocard_ramanujan_12th_theorem}
If $n!=a^2b=a^2(c^2+d),~d\le 2c,~n\in \mathbb{N}$ and $n!+1=(ac+e)^2$, then for the solution of Brocard-Ramanujan problem (c.f. (\ref{brocard_ramanujan_problem})), $1\leq d<1+\frac{2ce}{a}$.

{\bf Proof}: For natural number $n(>1)$, we know that 
\begin{eqnarray} 
	n!=a^2b=a^2(c^2+d),~d\le 2c.\label{brocard_ramanujan_relation_nacd}\\
\text{Then~}n!+1=a^2(c^2+d)+1=(ac+e)^2\label{brocard_ramanujan_relation_nace}\\
 \text{if~} a^2d+1=2ace+e^2.\label{brocard_ramanujan_relation_acde}\\
\text{For~} e=1, ad=2c.	\label{brocard_ramanujan_relation_acd_e=1}\\
\text{If~}e>1, a(ad-2ce)=e^2-1. \label{brocard_ramanujan_relation_acde_e>1}\\
\text{Thus~}a/(e^2-1)~\text{and~}ad-2ce>0~\text{i.e.~}e<\frac{ad}{2c}\label{brocard_ramanujan_relation_ad_2ce}\\
\text{Or,~} e<a~\text{as~}d\le 2c.\label{brocard_ramanujan_relation_e<a}
\end{eqnarray}	 
 We observe that $ad-2ce\ge a$ or $ad-2ce<a$. 
\begin{eqnarray} 
	\text{If~} ad-2ce\geq a~{i.e.~}a(d-1)\geq 2ce,\label{brocard_ramanujan_relation_ad-2ce>=a}\\
	\text{then~}d>1~\text{as~}ce>0\text{~and~}d>1+\frac{2ce}{a}\label{brocard_ramanujan_relation_d>1+2ce/a}
\end{eqnarray}
But then in (\ref{brocard_ramanujan_relation_acde_e>1}), L.H.S.$>a^2$ and R.H.S.$<a^2$. Thus (\ref{brocard_ramanujan_relation_acde}) and hence (\ref{brocard_ramanujan_relation_nace}) has no solution.
\begin{eqnarray} 
\text{If~} ad-2ce<a,~{then~}a(d-1)<2ce,\label{brocard_ramanujan_relation_ad-2ce<a}
\end{eqnarray}
When $d<1$ (i.e. $d=0$), then (\ref{brocard_ramanujan_relation_acde_e>1}) is always true, but $d=0$ contradicts (\ref{brocard_ramanujan_relation_ad_2ce}) as $ce>0$. Also $n!\neq a^2c^2$ if $n>1$. 

If $d=1$, $ce>0$ from (\ref{brocard_ramanujan_relation_ad-2ce<a}). Using (\ref{brocard_ramanujan_relation_ad_2ce}) and (\ref{brocard_ramanujan_relation_ad-2ce<a}), we get $0<a-2ce<a$. Thus (\ref{brocard_ramanujan_relation_acde}) and hence (\ref{brocard_ramanujan_relation_nace}) has solution if it satisfies (\ref{brocard_ramanujan_relation_acde_e>1}). From (\ref{brocard_ramanujan_relation_acde_e>1}), we have $1<d<1+\frac{2ce}{a}$ and $0<ad-2ce<a$. This means (\ref{brocard_ramanujan_relation_acde}) and hence (\ref{brocard_ramanujan_relation_nace}) has solution if it satisfies (\ref{brocard_ramanujan_relation_acde_e>1}). 

\subsubsection{Remark}
\label{brocard_ramanujan_5th_remark}
For $n\ge 6$; $a^2=a_1^2n^2$ if $n$ is not prime number. $a^2=a_1^2(n+1)^2$ if $n$ is a prime number ($n\neq 7$) but $\frac{n+1}{2}$ is not a prime number and $a^2=a_1^2(n+2)^2$ if $n$, $\frac{n+1}{2}$ are both prime numbers following the examples. If $n=7$, then $a^2=a_1^2(n+5)^2=1^212^2$. Moreover, $\sqrt{n!}=a(c+\epsilon_3),~2c\epsilon_3+\epsilon_3^2=d,~(0<\epsilon_3<1)$. Thus satisfying the Theorem \ref{brocard_ramanujan_3rd_theorem}, we have $k=ac+[a\epsilon_3],~\epsilon=a\epsilon_3-[a\epsilon_3]$.

1. $1!=1=1^2,~a=1,~b=1,~c=1,~d=0,~\epsilon_3=0,~k=1,~\epsilon=0 $.

2. $2!=2=1^2(1^2+1),~a=1,~b=1,~c=1,~d=1,~\epsilon_3\approx 0.414213562,~k=1,~\epsilon\approx 0.414213562$.

3. $3!=2\times 3=1^2(2^2+2),~a=1,~b=6,~c=2,~d=2,~\epsilon_3\approx 0.449489743,~k=2,~\epsilon\approx 0.449489743$.

4. $4!=2^2\times 2\times 3=2^2(2^2+2),~a=2,~b=6,~c=2,~d=2,~\epsilon_3\approx 0.449489743,~k=4,~\epsilon\approx 0.898979486$.

5. $5!=2^2\times 2\times 3\times 5=2^2(5^2+5),~a=2,~b=30,~c=5,~d=5,~\epsilon_3\approx 0.477225575,~k=10,~\epsilon\approx 0.95445115$.

6. $6!=4^2 3^2 5=2^2 6^25=12^2(2^2+1),~a=12,~a_1=2,~b=5,~c=2,~d=1,~\epsilon_3\approx 0.236067977,~k=26,~\epsilon\approx 0.83281573$.

7. $7!=4^2 3^2 5\times 7=1^2 12^2(5^2+10),~a=12,~a_1=1,~b=35,~c=5,~d=10,~\epsilon_3\approx 0.916079783,~k=70,~\epsilon\approx 0.992957397$.

8. $8!=8^2 3^2 2\times5\times 7=3^2 8^2(8^2+6),~a=24,~a_1=3,~b=70,~c=8,~d=6,~\epsilon_3\approx 0.366600265,~k=200,~\epsilon\approx 0.798406368$.

9. $9!=8^2 9^2 2\times5\times 7=8^2 9^2(8^2+6),~a=72,~a_1=8,~b=70,~c=8,~d=6,~\epsilon_3\approx 0.366600265,~k=602,~\epsilon\approx 0.39521908$.

10. $10!=16^2 9^2 5^2 7=72^2 10^2(2^2+3),~a=720,~a_1=72,~b=7,~c=2,~d=3,~\epsilon_3\approx 0.645751311,~k=1904,~\epsilon\approx 0.940943967$.

11. $11!=16^2 9^2 5^2 7\times 11=60^2 12^2(8^2+13),~a=720,~a_1=60,~b=77,~c=8,~d=13,~\epsilon_3\approx 0.774964387,~k=6317,~\epsilon\approx 0.974358922$.

12. $12!=32^2 9^2 5^2 3\times7\times 11=120^2 12^2(8^2+13),~a=1440,~a_1=120,~b=231,~c=15,~d=6,~\epsilon_3\approx 0.198684154,~k=21886,~\epsilon\approx 0.105181142$.

13. $13!=32^2 9^2 5^2 3\times7\times 11\times 13=96^2 15^2(8^2+13),~a=1440,~a_1=96,~b=3003,~c=54,~d=87,~\epsilon_3\approx 0.799635035,~k=83517,~\epsilon\approx 0.4744504$.

14. $14!=32^2 9^2 5^2 7^2 2\times 3\times 11\times 13=720^2 14^2(29^2+17),~a=10080,~a_1=720,~b=858,~c=29,~d=17,~\epsilon_3\approx 0.291637032,~k=295259,~\epsilon\approx 0.701280076$.

15. $15!=32^2 27^2 5^2 7^2 2\times 5\times 11\times 13=2016^2 15^2(37^2+61),~a=30240,~a_1=2016,~b=1430,~c=37,~d=61,~\epsilon_3\approx 0.815340802,~k=1143535,~\epsilon\approx 0.905863913$.

16. $16!=128^2 27^2 5^2 7^2 2\times 5\times 11\times 13=7560^2 16^2(37^2+61),~a=120960,~a_1=7560,~b=1430,~c=37,~d=61,~\epsilon_3\approx 0.815340802,~k=4574143,~\epsilon\approx 0.623455652$.

17. $17!=128^2 27^2 5^2 7^2 2\times 5\times 11\times 13\times 17=6720^2 18^2(155^2+285),~a=120960,~a_1=6720,~b=24310,~c=155,~d=285,~\epsilon_3\approx 0.815340802,~k=18859677,~\epsilon\approx 0.306253148$.

18. $18!=256^2 81^2 5^2 7^2 5\times 11\times 13\times 17=40320^2 18^2(110^2+55),~a=725760,~a_1=40320,~b=12155,~c=110,~d=55,~\epsilon_3\approx 0.249716553,~k=80014834,~\epsilon\approx 0.28550528$.

19. $19!=256^2 81^2 5^2 7^2 5\times 11\times 13\times 17\times 19=36288^2 20^2(480^2+545),~a=725760,~a_1=36288,~b=230945,~c=480,~d=545,~\epsilon_3\approx 0.567373008,~k=348776576,~\epsilon\approx 0.63428608$.

20. $20!=512^2 81^2 25^2 7^2 11\times 13\times 17\times 19=362880^2 20^2(214^2+393),~a=7257600,~a_1=362880,~b=46189,~c=214,~d=393,~\epsilon_3\approx 0.916262763,~k=1559776268,~\epsilon\approx 0.6287488$.

21. $21!=512^2 81^2 25^2 7^2 3\times 7\times 11\times 13\times 17\times 19=345600^2 21^2(984^2+1713),~a=7257600,~a_1=345600,~b=969969,~c=984,~d=1713,~\epsilon_3\approx 0.870042188,~k=7147792818,~\epsilon\approx 0.1836288$.

22. $22!=512^2 81^2 25^2 7^2 11^2 2\times 3\times 7\times 13\times 17\times 19=3628800^2 22^2(419^2+797),~a=79833600,~a_1=3628800,~b=176358,~c=419,~d=797,~\epsilon_3\approx 0.949997023,~k=33526120082,~\epsilon\approx 0.3353728$.

23. $23!=512^2 81^2 25^2 7^2 11^2 2\times 3\times 7\times 13\times 17\times 19\times 23=3326400^2 24^2(2014^2+38),~a=79833600,~a_1=3326400,~b=4056234,~c=2014,~d=38,~\epsilon_3\approx 0.00943394,~k=160785623545,~\epsilon\approx 0.392384$.

24. $24!=2048^2 243^2 25^2 7^2 11^2 7\times 13\times 17\times 19\times 23=79833600^2 24^2(822^2+355),~a=1916006400,~a_1=79833600,~b=676039,~c=822,~d=355,~\epsilon_3\approx 0.2159083844,~k=,~\epsilon\approx 0.5576576$.

\subsubsection{Lemma}
\label{brocard_ramanujan_2nd_lemma} 
For a given natural number $n>1$, we know $n!=a^2b=a^2(c^2+d),~d\leq 2c,~c<an$ and $n+1=c_1^2+d_1,~d_1\leq 2c_1$. Then $(n+1)!<a^2(3c^2c_1^2+2cc_1),~c_1<\sqrt{n}$.

{\bf Proof}: Given that $n!=a^2b=a^2(c^2+d),~d\leq 2c,~c<an$ and $n+1=c_1^2+d_1,~d_1\leq 2c_1$. Then for $d=2c,~d_1=2c_1;~(n+1)!=a^2(c^2+2c)(c_1^2+2c_1)=a^2(c^2c_1^2+2cc_1^2+2c_1c^2+4cc_1)=a^2(c^2c_1^2+2cc_1(c+c_1+2))<a^2(3c^2c_1^2+2cc_1)$ as $1+c+c_1<c+c_1$ when $cc_1-c-c_1+1>2$ i.e. $(c-1)(c_1-1)>2$. Again $n+1=c^2+2c_1$ gives $c_1<\sqrt{n}$ and $\sqrt{3}~cc_1<\sqrt{3}\sqrt{n}~an<n^2a$.

\subsection{Theorem}
\label{brocard_ramanujan_13th_theorem}
If $n!=a^2b=a^2(c^2+d),~d\le 2c,~n\in \mathbb{N}$ and $\sqrt{n!}=k+\epsilon$, then $2k\epsilon+\epsilon^2=a^2d-[a\epsilon_3](2k-[a\epsilon_3])=a^2d-[a\epsilon_3](2ac+[a\epsilon_3])=n!-k^2$, where $\sqrt{c^2+d}=c+\epsilon_3$. For the solution of Brocard-Ramanujan problem (c.f. (\ref{brocard_ramanujan_problem})), $2k\epsilon+\epsilon^2=a^2d-[a\epsilon_3](2k-[a\epsilon_3])=a^2d-[a\epsilon_3](2ac+[a\epsilon_3])=2k$ and $k=\frac{a^2d+[a\epsilon_3]^2}{2(1+[a\epsilon_3])}$.

{\bf Proof}: For natural number $n(>1)$, we know that 
\begin{eqnarray} 
	n!=a^2(c^2+d),~d\le 2c.\label{brocard_ramanujan_relation_nacd2}~~~~~~~~~\\
	\text{Then~}\sqrt{n!}=a(c+\epsilon_3)~(\text{where~}2c\epsilon_3+\epsilon_3^2=d)\label{brocard_ramanujan_relation_nac_epsilon3}~~~~~~~~~\\
	=(ac+[a\epsilon_3])+(a\epsilon_3-[a\epsilon_3])=k+\epsilon,~k=ac+[a\epsilon_3],~\epsilon=a\epsilon_3-[a\epsilon_3].\label{brocard_ramanujan_relation_nac_epsilon3_k_epsilon}~~~~~~~~~\\
	\text{Now,~} 2k\epsilon+\epsilon^2=2(ac+[a\epsilon_3])(a\epsilon_3-[a\epsilon_3])+(a\epsilon_3-[a\epsilon_3])^2=a\epsilon_3(2ac+a\epsilon_3)-[a\epsilon_3](2ac+[a\epsilon_3])~~~~~~\nonumber\\=a^2(2c\epsilon_3+\epsilon_3^2)-[a\epsilon_3](2ac+[a\epsilon_3])~~~~~~~~~\nonumber\\=a^2d-[a\epsilon_3](2ac+[a\epsilon_3])
	=a^2d-[a\epsilon_3](2k-[a\epsilon_3])~~~~~~~~~ \label{brocard_ramanujan_relation_2kepsilon_epsiloin^2}\\=a^2d-(k-ac)(k+ac)=a^2d-k^2+a^2c^2=a^2(c^2+d)-k^2=n!-k^2.~~~~~~~~~\label{brocard_ramanujan_relation_2kepsilon_epsiloin^2=n!-k^2}
	\end{eqnarray}	 
Moreover, we know that $n!=k^2+2k\epsilon+\epsilon^2$ i.e. $2k\epsilon+\epsilon^2=n!-k^2$. Thus for the solution of Brocard-Ramanujan problem, we have from (\ref{brocard_ramanujan_relation_2kepsilon_epsiloin^2}) 
\begin{eqnarray} 
2k\epsilon+\epsilon^2=a^2d-[a\epsilon_3](2ac+[a\epsilon_3])=a^2d-[a\epsilon_3](2k-[a\epsilon_3])=2k.\label{brocard_ramanujan_relation_2kepsilon_epsilon^2=2k2nd}\\
	\text{Then~}k=ac+[a\epsilon_3]=\frac{a^2d+[a\epsilon_3]^2}{2(1+[a\epsilon_3])}\label{brocard_ramanujan_relation_akdepsilon3}
\end{eqnarray}
For the solution of Brocard-Ramanujan problem, we have the following examples:

1. $4!=2^26=2^2(2^2+2),~a=2,~c=2,~d=2,~~\epsilon_3\approx 0.449489743,~k=4,~\epsilon\approx 0.898979486$.

2. $5!=2^2\times 2\times 3\times 5=2^2(5^2+5),~a=2,~b=30,~c=5,~d=5,~\epsilon_3\approx 0.477225575,~k=10,~\epsilon\approx 0.95445115$.

3. $7!=4^2 3^2 5\times 7=1^2 12^2(5^2+10),~a=12,~b=35,~c=5,~d=10,~\epsilon_3\approx 0.916079783,~k=70,~\epsilon\approx 0.992957397$.
\subsubsection{Remark}
\label{brocard_ramanujan_6th_remark}
We know that $k$ is strictly increasing (if $n>2$) and $a$ is increasing with the increase value of $n$. Thus from (\ref{brocard_ramanujan_relation_akdepsilon3}), $k$ will be increasing if $a,~d$ (and hence may be $[a\epsilon_3]$) are increasing. We can observe that $d=c$ for first and second solution of Brocard-Ramanujan problem. For the third solution, it is $d=2c$. Also, there is no other solution for $n\leq 10^9$. Moreover, $n!=k^2+2k$ for the solution and $n!<k^2+2k$ for non-solution of (\ref{brocard_ramanujan_problem}). Thus $a,~d,~\epsilon_3$ are increasing for maximum increase of $k$ to the solution, however maximum value of $d$ is $2c$. Hence, we may claim here also that there will be no other solution of  (\ref{brocard_ramanujan_problem}).

\subsubsection{Corollary}
\label{brocard_ramanujan_4th_corollary}
$n!=k^2+a^2d-[a\epsilon_3](2ac+[a\epsilon_3])=k^2-[a\epsilon_3](2k-[a\epsilon_3])=k^2+(n!-k^2).$

\subsection{Theorem}
\label{brocard_ramanujan_14th_theorem}
If $n!=a^2b=a^2(c^2+d),~d\le 2c,~n\in \mathbb{N}$, then $n!+1=(ac+e)^2$ if $a^2d+1=e(2ac+e)$ for the solution of Brocard-Ramanujan problem and $(ac+e)^2-1<n!+1<(ac+e)^2$ for the non-solution of (\ref{brocard_ramanujan_problem}). By division algorithm $a^2d+1=q\times 2ac+r,~r<2ac$. Then (\ref{brocard_ramanujan_problem}) has solution if and only if $a^2c^2+2acq+r=\text{square}$. Moreover, the value of $e$ is $e=q-[x]$ where $x^2-2(ac+q)x+(q^2-r)=0$. 

{\bf Proof}: We have  
\begin{eqnarray} 
	a^2d+1=q\times 2ac+r,~(r<2ac).\label{brocard_ramanujan_theo14_eq1}~~~~~~~~~\\=(q-x)\{2ac+(q-x)\}+r+2acx-(q-x)^2.\label{brocard_ramanujan_theo14_eq2}~~~~~~~~~\\
	\text{Then~}n!+1=(ac+q-x)^2~\text{if~}2acx-(q-x)^2+r=0.\label{brocard_ramanujan_theo14_eq3}~~~~~~~~~\\
	\text{Or,~}x^2-2(ac+q)x+(q^2-r)=0.\label{brocard_ramanujan_theo14_eq4}~~~~~~~~~\\
	\text{Then,~} x=\frac{2(ac+q)\pm\sqrt{4(ac+q)^2-4(q^2-r)}}{2}=(ac+q)\pm \sqrt{a^2c^2+2acq+r}.\label{brocard_ramanujan_theo14_eq5}~~~~~~~~~
\end{eqnarray}	 
For the solution of (\ref{brocard_ramanujan_problem}), we need $x\in\mathbb{N}$.
\begin{eqnarray} 
\text{Then,~}a^2c^2+2acq+r=\text{square}=f^2,~f\in\mathbb{N}	;\label{brocard_ramanujan_theo14_eq6}~~~~~~~~~\\n!+1=(ac+e)^2=f^2=m^2,~e=q-x,~f=ac+(q-x).\label{brocard_ramanujan_theo14_eq7}~~~~~~~~~
\end{eqnarray}

For the non-solution of (\ref{brocard_ramanujan_problem}), we have $x\notin\mathbb{N}$. Then $x=(ac+q)-\sqrt{a^2c^2+2acq+r}=x_1+x_2,~0<x_2<1$; $(ac+e)^2-1<n!+1<(ac+e)^2$, $e=q-x_1=q-[x]$. Moreover, $a^2d+1=(q-x)\{2ac+(q-x)\}~\text{as~}r+2acx-(q-x)^2=0$. So $a^2d+1=(q-x_1-x_2)\{2ac+(q-x_1-x_2)\}=(q-x_1)\{2ac+(q-x_1)\}-2acx_2-2x_2(q-x_1)+x_2^2$. Thus $acx_2+2x_2(q-x_1)-x_2^2<2(ac+e)-1$. Some examples of this theorem have been presented below.

1. $2!=1^2(1^1+1),~a=1,~c=1,~d=1,~a^2d+1=2=1\times 2+0,~x=2-\sqrt{3},~e=q-[x]=1,~ac+e=2,~1^2<2!+1<2^2$.

2. $3!=1^2(2^2+2),~a=1,~c=2,~d=2,~a^2d+1=3=0\times 4+3,~x=2-\sqrt{7},~e=q-[x]=1,~ac+e=3,~2^2<3!+1<3^2$.

3. $4!=2^2(2^2+2),~a=2,~c=2,~d=2,~a^2d+1=9=1\times 8+1,~x=5-\sqrt{5^2},~e=q-x=1,~ac+e=5,~4!+1=5^2$.

4. $5!=2^2(5^2+5),~a=2,~c=5,~d=5,~a^2d+1=21=1\times 20+1,~x=11-\sqrt{11^2},~e=q-x=1,~ac+e=11,~5!+1=11^2$.

5. $6!=12^2(2^2+1),~a=12,~c=2,~d=1,~a^2d+1=145=3\times 48+1,~x=27-\sqrt{27^2-8},~e=q-[x]=3,~ac+e=27,~26^2<6!+1<27^2$.

6. $7!=12^2(5^2+10),~a=12,~c=5,~d=10,~a^2d+1=1441=12\times 120+1,~x=72-\sqrt{71^2},~e=q-x=11,~ac+e=71,~7!+1=71^2$.

7. $8!=24^2(8^2+6),~a=24,~c=8,~d=6,~a^2d+1=3457=9\times 384+1,~x=201-\sqrt{200^2+321},~e=q-[x]=9,~ac+e=201,~200^2<8!+1<201^2$.

8. $17!=120960^2(155^2+285),~a=120960,~c=155,~d=285,~a^2d+1=3=111205\times 37497600+6048001,~x=18860005-\sqrt{120960^2155^2+2\times 120960\times155\times 6048001}\approx 327.693746824,~e=q-[x]=110878,~ac+e=18859678,~18859677^2<17!+1<18859678^2$.

\subsubsection{Corollary}
\label{brocard_ramanujan_5th_corollary}
If $d=2c$, then $a^2c^2+2a^2c+1=f^2,~f^2-1=a^2(f_1^2-1);~f,~f_1\in\mathbb{N}$ for the solution of (\ref{brocard_ramanujan_problem}).

{\bf Proof}: 
\begin{eqnarray} 
\text{For~}d=2c,~a^2d+1=2a^2c+1=a\times 2ac+1,~q=a,~r=1.\label{brocard_ramanujan_cor5_eq1}~~~~~~~~~\\x=(ac+a)-\sqrt{a^2c^2+2a^2c+1}.\label{brocard_ramanujan__cor5_eq2}~~~~~~~~~
\end{eqnarray}	 
For the solution of (\ref{brocard_ramanujan_problem}), we need $x\in\mathbb{N}$.
\begin{eqnarray} 
	\text{Then~}a^2c^2+2a^2c+1=f^2,~f\in\mathbb{N}.\label{brocard_ramanujan_cor5_eq3}~~~~~~~~~~~~~~~~~~~~~~~~\\\text{Or,~}c=\frac{-2a^2+\sqrt{4a^4+4a^2(f^2-1)}}{2a^2}=-1+\frac{\sqrt{a^2+f^2-1}}{a}.\label{brocard_ramanujan__cor5_eq4}~~~~~~~~~
\end{eqnarray}
Thus for the solution of (\ref{brocard_ramanujan_problem}), $a^2$ divides $f^2-1$, 
\begin{eqnarray} 
f^2-1=a^2(f_1^2-1),~f_1\in\mathbb{N}\label{brocard_ramanujan_cor5_eq5}~~~~~~~~~\\\text{and~}x=(ac+a)-f,~c=-1+f_1.\label{brocard_ramanujan__cor5_eq6}~~~~~~~~~
\end{eqnarray}
When $n=7$, $a=12,~7!=12^235=12^2(5^2+2\times 5)=12^2(6^2-1)=71^2-1,~f=71,~f_1=6$. Moreover, $5040=71^2-1=70\times 72=12^2(6^2-1)=60\times 84$. Thus $7!+1=71^2$.

We also can find that $f=97,~f^2-1=96\times 98=14^2(7^2-1)=14^2(6^2+2\times 6)$. But there does not exist $n\in\mathbb{N}$ such that $n!=14^2(6^2+12)$. We can generate infinite number of such solutions of $f^2-1=a^2(f_1^2-1)$ which follow unique pattern as given below.

$4^2(1^2+2\times 1)=4^2(2^2-1)=7^2-1$,

$6^2(2^2+2\times 2)=6^2(3^2-1)=17^2-1$,

$8^2(3^2+2\times 3)=8^2(4^2-1)=31^2-1$,

$10^2(4^2+2\times 4)=10^2(5^2-1)=49^2-1$

$16^2(7^2+2\times 7)=16^2(8^2-1)=127^2-1$,

$18^2(8^2+2\times 8)=18^2(9^2-1)=161^2-1$,

$20^2(9^2+2\times 9)=20^2(10^2-1)=199^2-1$.

Thus we have the unique pattern as
\begin{eqnarray} 
	(2c+2)^2(c^2+2c)=(2c+2)^2\{(c+1)^2-1\}=\{2c(c+2)+1\}^2-1,~c\in\mathbb{N}.\label{brocard_ramanujan_cor5_eq7}~~~~~~~~~
\end{eqnarray}
Also from (\ref{brocard_ramanujan__cor5_eq2}), (\ref{brocard_ramanujan_cor5_eq3}) and (\ref{brocard_ramanujan_cor5_eq7}) using $a=2c+2$, we can obtain 
\begin{eqnarray} 
x=(c+1)(2c+2)-\sqrt{2^2c^2(c+2)^2+4c(c+2)+1}=2(c+1)^2-2c(c+2)=1.\label{brocard_ramanujan_cor5_eq8}~~~~~~~~~
\end{eqnarray}
Analyzing the section \ref{brocard_ramanujan_5th_remark} Remark, it can be discovered that $a$ is even for $n>3$, $2c+2<a$ for $n>7$ and $nc<<a$ for large $n$ ($n\geq 25$). Hence, there will be no other solution of (\ref{brocard_ramanujan_problem}) following pattern of equation (\ref{brocard_ramanujan_cor5_eq7}) except the only solution $7!+1=71^2$ when $d=2c$. However, $a\leq c$ when $n\leq 5$. Then for the solution of (\ref{brocard_ramanujan_problem}), we need $d<2c$. We have already observed that there are two solutions $4!+1=5^2$ and $5!+1=11^2$ of (\ref{brocard_ramanujan_problem}) when $d=c$. 

\subsubsection{Corollary}
\label{brocard_ramanujan_6th_corollary}
If $d=c$, then $a^2c^2+a^2c+1=f^2,~4(f^2-1)=a^2(f_1^2-1);~f,~f_1\in\mathbb{N}$ for the solution of (\ref{brocard_ramanujan_problem}).

{\bf Proof}: 
\begin{eqnarray} 
	\text{For~}d=c,~a^2d+1=a^2c+1=\frac{a}{2}\times 2ac+1,~q=\frac{a}{2},~r=1.\label{brocard_ramanujan_cor6_eq1}~~~~~~~~~\\x=(ac+a)-\sqrt{a^2c^2+2a^2c+1}.\label{brocard_ramanujan__cor6_eq2}~~~~~~~~~
\end{eqnarray}	 
For the solution of (\ref{brocard_ramanujan_problem}), we need $x\in\mathbb{N}$.
\begin{eqnarray} 
	\text{Then~}a^2c^2+a^2c+1=f^2,~f\in\mathbb{N}.\label{brocard_ramanujan_cor6_eq3}~~~~~~~~~~~~~~~~~~~~~~~~~~\\\text{Or,~}c=\frac{-a^2+\sqrt{a^4+4a^2(f^2-1)}}{2a^2}=-\frac{1}{2}+\frac{\sqrt{a^2+4(f^2-1)}}{2a}.\label{brocard_ramanujan__cor6_eq4}~~~~~~~~~
\end{eqnarray}
Thus for the solution of (\ref{brocard_ramanujan_problem}), $a^2$ divides $4(f^2-1)$, 
\begin{eqnarray} 
	4(f^2-1)=a^2(f_1^2-1),~f_1\in\mathbb{N}\label{brocard_ramanujan_cor6_eq5}~~~~~~~~~\\\text{and~}x=\left(ac+\frac{a}{2}\right)-f,~c=-\frac{1}{2}+\frac{f_1}{2},~f_1~\text{is~odd}.\label{brocard_ramanujan__cor6_eq6}~~~~~~~~~
\end{eqnarray}
When $n=4$, $a=2,~4!=2^26=2^2(2^2+2)=5^2-1,~2^2(f^2-1)=a^2(f_1^2-1),~f=f_1=5$; and for $n=5$, $a=2,~5!=2^230=2^2(5^2+5)=11^2-1,~2^2(f^2-1)=a^2(f_1^2-1),~f=f_1=11$. 

We also can find that $2^2(6^2+6)=13^2-1,~a=2,~f=f_1=13$. But there does not exist $n\in\mathbb{N}$ such that $n!=2^2(6^2+6)$. We can generate infinite number of such solutions of $2^2(f^2-1)=a^2(f_1^2-1)$ which follow unique pattern as given below.

$2^2(1^2+1)=3^2-1,~a=2,~f=f_1=3$,

$2^2(3^2+3)=7^2-1,~a=2,~f=f_1=7$,

$2^2(4^2+4)=7^2-1,~a=2,~f=f_1=9$,

$2^2(7^2+7)=15^2-1,~a=2,~f=f_1=15$,

$2^2(8^2+8)=17^2-1,~a=2,~f=f_1=17$,

$2^2(9^2+9)=19^2-1,~a=2,~f=f_1=19$.

Thus we have the unique pattern as
\begin{eqnarray} 
	2^2(c^2+c)=(4c^2+4c)=(2c+1)^2-1,a=2,~f=f_1=2c+1,~c\in\mathbb{N}.\label{brocard_ramanujan_cor6_eq7}~~~~~~~~~
\end{eqnarray}
Analyzing the section \ref{brocard_ramanujan_5th_remark} Remark, it can be discovered that $a$ is even for $n>3$, $a>2c+2$ for $n>7$ and $a>>nc$ for large $n$ ($n>50$). Hence, there will be no other solution of (\ref{brocard_ramanujan_problem}) following pattern of equation (\ref{brocard_ramanujan_cor6_eq7}) except the two solutions $4!+1=5^2$ and $5!+1=11^2$ when $d=c$.

\subsubsection{Corollary}
\label{brocard_ramanujan_7th_corollary}
By further thoroughly investigation of the pattern of (\ref{brocard_ramanujan_cor5_eq7}), we can generate infinite number of similar patterns.

{\bf Proof}:
Following the further investigation of the pattern of (\ref{brocard_ramanujan_cor5_eq7}), we can get 
\begin{eqnarray} 
\{2(2c^2+4c)+2\}^2\{(2c^2+4c)^2+2(2c^2+4c)\}~~~~~~~~~~~~~~~~~~~~~~~~~~~~~\nonumber\\=\{8(c^2+2c)^2+8(c^2+4c)+2\}\{8(c^2+2c)^2+8(c^2+4c)\}~~~~~~~~~\nonumber\\ \text{Or,~}4^2\{2c(c+2)+1\}^2(c+1)^2(c^2+2c)\}=\{8c(c+1)^2(c+2)+1\}^2-1.\label{brocard_ramanujan_cor7_eq1}~~~~~~~~~
\end{eqnarray}
Here $a=4(c+1)\{2c(c+2)+1\},~d=2c,~q=a,~r=1$ and then 
\begin{eqnarray} 
x=(ac+q)-\sqrt{a^2c^2+2acq+r}\nonumber~~~~~~~~~~~~~~~~~~~~~~~~~~~~~~~~~~~~~~~~~~~~~~~~~~~~~~~~~~~~~~~\\=4(c+1)^2\{2c(c+2)+1\}-\sqrt{64c^8+512c^7+1664c^6+2816c^5+2640c^4+1344c^3+336c^2+32c+1}\nonumber\\=4(c-1)(c+3)+15.\label{brocard_ramanujan_cor7_eq2}~~~~~~~~~~~~~~~~~~~~~~~~~~~~~~~~~~~~~~~~~~~~~~~~~~~~~~~~~~~~~~~~~~~~~~~~~~~~~~~~~~~~~~~~~~~~~~~~~~~~
\end{eqnarray}

Following (\ref{brocard_ramanujan_cor7_eq1}), we again obtain 
\begin{eqnarray}  
	\{16c(c+1)^2(c+2)+2\}^2\{8c(c+1)^2(c+2)\}\{8c(c+1)^2(c+2)+2\}\nonumber~~~~~~~~~~~~~~~~~~~~~~~~~~~~~~~~~~~~\\= \{128c^2(c+1)^4(c+2)^2+32c(c+1)^2(c+2)+2\}\{128c^2(c+1)^4(c+2)^2+32c(c+1)^2(c+2)\}~~~~\label{brocard_ramanujan_cor7_eq3}
\end{eqnarray}
\begin{eqnarray} 
	\text{Or,~}8^2\{8c(c+1)^2(c+2)+1\}^2\{2c(c+2)+1\}^2(c+1)^2(c^2+2c)\nonumber~~~~~~~~~~~~~~~~~~~~~~~~~~~~~~~~~\\= \{32c(c+1)^2(c+2)\{2c(c+2)+1\}^2+1\}^2-1.\label{brocard_ramanujan_cor7_eq4}~~~~~~~~~~~~~~~~~~~~~~~~~~~~~~~~~~~~~~~~~~~	
\end{eqnarray}	
We can repeat this process.

\subsubsection{Corollary}
\label{brocard_ramanujan_8th_corollary}
Investigating the pattern of (\ref{brocard_ramanujan_cor5_eq7}), we can formulate infinite number of similar patterns for $a^2(c^2+d)$.

{\bf Proof}:
Following the pattern of (\ref{brocard_ramanujan_cor5_eq7}), we can evaluate
\begin{eqnarray} 
	\{2(2c+1)\}^2(4c^2+4c)=(8c^2+8c+1)(8c^2+8c).~~~~~~~~~\\ \text{Or,~}4^2(2c+1)^2(c^2+c)=\{8c(c+1)+1\}^2-1.\label{brocard_ramanujan_cor8_eq1}~~~~~~~~~
\end{eqnarray}
Here $a=4(2c+1),~d=c,~q=\frac{a}{2},~r=1$ and then 
\begin{eqnarray} 
	x=(ac+q)-\sqrt{a^2c^2+2acq+r}=2(2c+1)^2-\sqrt{64c^4+128c^3+80c^2+16c+1}=1. \label{brocard_ramanujan_cor8_eq2}~~
\end{eqnarray}

Following (\ref{brocard_ramanujan_cor8_eq1}), we again obtain 
\begin{eqnarray}  
	\{16c(c+1)^2(c+2)+2\}^2\{8c(c+1)^2(c+2)\}\{8c(c+1)^2(c+2)+2\}\nonumber~~~~~~~~~~~~~~~~~~~~~~~~~~~~~~~~~~~~\\= \{128c^2(c+1)^4(c+2)^2+32c(c+1)^2(c+2)+2\}\{128c^2(c+1)^4(c+2)^2+32c(c+1)^2(c+2)\}~~~~\label{brocard_ramanujan_cor8_eq3}
\end{eqnarray}
\begin{eqnarray} 
	\text{Or,~}8^2\{8c(c+1)^2(c+2)+1\}^2\{2c(c+2)+1\}^2(c+1)^2(c^2+2c)\nonumber~~~~~~~~~~~~~~~~~~~~~~~~~~~~~~~~~\\= \{32c(c+1)^2(c+2)\{2c(c+2)+1\}^2+1\}^2-1.\label{brocard_ramanujan_cor8_eq4}~~~~~~~~~~~~~~~~~~~~~~~~~~~~~~~~~~~~~~~~~~~	
\end{eqnarray}	
We can also repeat this process.

\subsubsection{Corollary}
\label{brocard_ramanujan_9th_corollary}
Investigating the pattern of (\ref{brocard_ramanujan_cor5_eq7}), we can formulate infinite number of similar patterns for $a^2(c^2+1)$.

{\bf Proof}:
Following the pattern of (\ref{brocard_ramanujan_cor5_eq7}), we have
\begin{eqnarray} 
	4c^2(c^2+1)=(2c^2+1)^2-1.\label{brocard_ramanujan_cor9_eq1}~~~~~~~~~
\end{eqnarray}
Again repeating this process, we can also obtain
\begin{eqnarray} 
16c^2(2c^2+1)^2(c^2+1)=\{8c^2(c^2+1)\}^2-1.	 \label{brocard_ramanujan_cor9_eq2}~~
\end{eqnarray}

\subsubsection{Corollary}
\label{brocard_ramanujan_10th_corollary}
Investigating the pattern of (\ref{brocard_ramanujan_cor5_eq7}), we can obtain infinite number of similar patterns for $a^2(c^2+2)$.

{\bf Proof}:
The pattern of (\ref{brocard_ramanujan_cor5_eq7}) helps us to formulate
\begin{eqnarray} 
	c^2(c^2+2)=(c^2+1)^2-1.\label{brocard_ramanujan_cor10_eq1}~~~~~~~~~
\end{eqnarray}
Repeating this process, we can also obtain
\begin{eqnarray} 
	4c^2(c^2+1)^2(c^2+2)=\{2c^2(c^2+2)+1\}^2-1.	 \label{brocard_ramanujan_cor10_eq2}~~
\end{eqnarray}

\subsubsection{Note}
 The patterns in sections \ref{brocard_ramanujan_7th_corollary}-\ref{brocard_ramanujan_10th_corollary} do not provide any new solution of Brocard-Ramanujan problem.

\subsection{Theorem}
\label{brocard_ramanujan_15th_theorem}
We can claim that Brocard-Ramanujan problem has no further solution using regular growth of $k$ and uncertain growth of $\epsilon$.

{\bf Proof}: By theorem (\ref{brocard_ramanujan_3rd_theorem}) and corollary (\ref{brocard_ramanujan_2nd_corollary}), we know that for the solution of Brocard-Ramanujan problem (\ref{brocard_ramanujan_problem}); $m-1=k=[\sqrt{n!}]$, $\sqrt{n!}=k+\epsilon$, where
$\epsilon=\sqrt{n!}-[\sqrt{n!}]$ i.e. $0<\epsilon <1$. Then for the solution of (\ref{brocard_ramanujan_problem}), the natural number $\epsilon (2k+\epsilon)$ must satisfy $\epsilon (2k+\epsilon)=2k$; $k=\frac{\epsilon^2}{2(1-\epsilon)}$ and
$\epsilon=\sqrt{k^2+2k}-k$ should be distinct and monotonic increasing for all the solutions of (\ref{brocard_ramanujan_problem}). 

Moreover, for any natural number $n$, $[\sqrt{n!}]=k$. Then for the natural number $n+1$, we have $[\sqrt{(n+1)!}]=[k\sqrt{n+1}]\geq k[\sqrt{n+1}]$. Thus for every growth of $n$, $k$ is growing with multiplicity $\sqrt{n+1}$. However, $\epsilon$ may increase or decrease as shown with examples in remark \ref{brocard_ramanujan_5th_remark}. Also, both the growths should satisfy $k=\frac{\epsilon^2}{2(1-\epsilon)}$ and
$\epsilon=\sqrt{k^2+2k}-k$.

According to Bruce Berndt and William Galway \cite{Berndt}, there is no new solution of the problem based on computations up to $n=10^9$. Thus for $n\geq 10^9$, we have huge growth of $k$ by regular growth with multiplicity $\sqrt{n+1}$, whereas growth of $\epsilon$ is uncertain such that its growth can not match the growth of $k$ to satisfy Brocard-Ramanujan problem (\ref{brocard_ramanujan_problem}). Hence, we can claim that Brocard-Ramanujan problem has no further solution.

\subsection{Conclusion}
The study explores multiple results on Brocard-Ramanujan
problem	(\ref{brocard_ramanujan_problem}). For natural number $n$, if $\sqrt{n!}=k+\epsilon,~n>1,~0<\epsilon<1$; then it has solution if and only if $n!=k(k+2)$ and $\epsilon,~k$  are strictly monotonic increasing. The problem has only finitely many solutions which is not based on any conjecture or previous study on the Brocard-Ramanujan problem. For the new solution when $n\ge 10^5$, the value of $\epsilon$ should be more than $0.999 \cdots 905915$ (digit 0 is coming after 228287 numbers of 9 digit taking more than 66 pages in LibreOffice Writer) indicating almost impossibility of new solution. For $n\geq 10^9$, I am unable to calculate the said numbers of 9 digit in the value of $\epsilon$ in my personal laptop (with 8GB Ram) using MATHEMATICA 8.

For the solution of (\ref{brocard_ramanujan_problem}), the natural number $\epsilon (2k+\epsilon)$ must satisfy $\epsilon (2k+\epsilon)=2k$ and conversely. Moreover, $k=\frac{\epsilon^2}{2(1-\epsilon)}$ and $\epsilon=\sqrt{k^2+2k}-k$ should be distinct and monotonic increasing for all the solutions of (\ref{brocard_ramanujan_problem}) and $n!<k(k+2)$ if $n$ is not a solution of Brocard-Ramanujan problem (c.f. equation
(\ref{brocard_ramanujan_problem})).

If $\sqrt{n!}=k+\epsilon$ for solution of Brocard-Ramanujan problem (c.f. equation (\ref{brocard_ramanujan_problem})), then $(\epsilon,k)$ lies on hyperbola.

If $n!=a^2b=a^2(c^2+d),~d\le 2c,~n\in \mathbb{N}$, $\sqrt{n!}=k+\epsilon$, $\sqrt{c^2+d}=c+\epsilon_3$; then for the solution of Brocard-Ramanujan problem, $k=\frac{a^2d+[a\epsilon_3]^2}{2(1+[a\epsilon_3])}$. Moreover, if $n!+1=(ac+e)^2$, then $a^2d+1=e(2ac+e)$ for the solution of Brocard-Ramanujan problem and $(ac+e)^2-1<n!+1<(ac+e)^2$ for the non-solution of (\ref{brocard_ramanujan_problem}). By division algorithm $a^2d+1=q\times 2ac+r,~r<2ac$ and (\ref{brocard_ramanujan_problem}) has solution if and only if $a^2c^2+2acq+r=\text{square}$. Also, the value of $e$ is $e=q-[x]$ where $x^2-2(ac+q)x+(q^2-r)=0$.

Finally, using the results of sections \ref{brocard_ramanujan_3rd_theorem}, \ref{brocard_ramanujan_3rd_remark}, \ref{brocard_ramanujan_10th_theorem}, \ref{brocard_ramanujan_6th_remark}, \ref{brocard_ramanujan_13th_theorem}, \ref{brocard_ramanujan_14th_theorem},\ref{brocard_ramanujan_5th_corollary},\ref{brocard_ramanujan_6th_corollary} \ref{brocard_ramanujan_15th_theorem}; it has been claimed  to discover that the problem has no further solution.
\vspace{1cm}

{\bf Acknowledgment:} {\it I am grateful to the University Grants Commission (UGC), New Delhi for awarding the Dr. D. S. Kothari Post Doctoral Fellowship from 9th July, 2012 to 8th July, 2015 at Indian Institute of Technology (BHU), Varanasi. The self training for this type of investigations was continued during the period. I am also acknowledge gratitude to my wife Dr. Ankita Chaturvedi who inspired me to complete this study. The problem was initiated at The LNM Institute of Information Technology, Jaipur, Rajasthan India.}

\end{document}